\documentstyle{amsppt} 
\magnification 1200
\overfullrule=0pt
\define\id{\operatorname{id}} 
 \define\gr{\operatorname{gr}} 
\define\co{\operatorname{co}} 
\define\ord{\operatorname{ord}}  
\define\End{\operatorname{End}} 
\define\length{\operatorname{length}} 
\define\Alg{\operatorname{Alg}}  \define\Ss{\operatorname{\Cal
S}} \define\varep{\varepsilon} \define\hq{\bold
h(q,m)}

\define\Aq{{\Cal A}\left(\Gamma, \Cal R, \Cal D\right)}
\define\bq{{\Cal B}\left(M, N, q, \lambda\right)}
\define\bqt{{\Cal B}\left(M, N, q, \widetilde{\lambda}\right)}

\topmatter  
\title
\nofrills {Lifting of Quantum Linear Spaces and Pointed
Hopf Algebras of order $ p^3$} \endtitle  
\rightheadtext{Lifting of Quantum Linear Spaces}
\author  {Nicol\'as
Andruskiewitsch and Hans-J\"urgen Schneider} \endauthor \address N.
Andruskiewitsch:  FAMAF. UNC. (5000) C. Universitaria. C\'ordoba.  Argentina
\endaddress \email  andrus\@mate.uncor.edu \endemail  \address H.-J.
Schneider: Mathematisches Seminar der Universit\"at  M\"unchen.
Theressienstr. 39. (80333) M\"unchen. Germany \endaddress
\email hanssch\@rz.mathematik.uni-muenchen.de\endemail
\thanks{(N.A.) Forschungsstipendiat der Alexander von
  Humboldt-Stiftung. Also partially supported by CONICET,
CONICOR, SeCYT (UNC)  and
  FaMAF (Rep\'ublica Argentina)m and by TWAS (Trieste)} \endthanks
  \endkeywords\abstract We propose the following
principle to study pointed Hopf algebras, or more generally, Hopf algebras
whose coradical is a Hopf subalgebra. Given such a Hopf algebra $A$,
consider its coradical filtration and the associated graded coalgebra $\gr
A$. Then $\gr A$ is a graded Hopf algebra, since the coradical $A_{0}$ of $A$
is a Hopf subalgebra. In addition, there is a projection $\pi: \gr A\to A_{0}$;
let $R$ be the algebra of coinvariants of $\pi$. Then, by a result of Radford
and Majid, $R$ is a braided Hopf algebra and $\gr A$ is the bosonization (or
biproduct) of $R$ and $A_{0}$: $\gr A \simeq  R\# A_{0}$. The principle we
propose to study $A$ is first to study $R$, then to transfer the information
to $\gr A$ via bosonization, and finally to lift to $A$. In this article, we
apply this principle to the situation when $R$ is the simplest braided Hopf
algebra: a quantum linear space. As consequences of our technique, we
obtain the classification of pointed Hopf algebras of order $p^{3}$ ($p$ an
odd prime) over an algebraically closed field of characteristic zero; with
the same hypothesis, the characterization of the  pointed Hopf algebras
whose coradical is abelian  and has index $p$ or $p^{2}$; and an infinite
family of pointed, non-isomorphic, Hopf algebras of the same dimension.
This last result gives a negative answer to a conjecture of I. Kaplansky.
\endabstract \endtopmatter 

\document  

\subhead \S 0. Introduction\endsubhead 
We assume for simplicity of the exposition that our ground-field $k$ is
algebraically closed of  characteristic 0; many results below are valid under
weaker hypotheses. Let $A$ be a non-cosemisimple  Hopf algebra whose
coradical $A_{0}$ is a Hopf subalgebra; for instance, $A$ is pointed, that is
all simple subcoalgebras are one dimensional.  Let   $$A_{0} \subseteq A_{1}
\subseteq \dots \subseteq A$$ be the coradical filtration of $A$, see
\cite{M, Chapter 5}. This is a coalgebra filtration and we consider the
associated graded coalgebra $\gr A = \oplus_{n \ge 0} \gr A(n)$,  $\gr A(n) =
A_{n}/ A_{n-1}$, where $A_{-1} = 0$. Since $A_{0}$ is a Hopf
subalgebra, $\gr A$ is a graded Hopf algebra and the zero term of
its own coradical filtration is $\gr A(0) = A_{0}$, which is  a Hopf
subalgebra of $\gr A$.  Let us denote $B = \gr A$, $H = \gr A(0)$. 
Let $\gamma: H \to B$ be the inclusion and let $\pi: B\to H$
be the projection with kernel $\oplus_{n \ge 1} \gr A(n)$. Then
$\pi$ is a Hopf algebra retraction of $\gamma$. We can describe
the situation in the following way:
 $$  R \hookrightarrow B  \overset \pi  \to{\underset \gamma
\to\rightleftarrows} H, \tag 0.1$$ 
where $R = B^{\co H} = \{a \in B: (\id \otimes \pi)
\Delta(a) = a\otimes 1\}$.  The setting (0.1) was first considered 
by Radford \cite{R3}; Majid  presented it in categorical terms 
\cite{Mj}. It turns out that $R$ is a  Hopf algebra  in the braided 
category  ${}_H^H\Cal{YD}$ of Yetter-Drinfeld modules over $H$;
we shall say "braided Hopf  algebra", for short. Moreover, $B$ can
be recovered as the biproduct (or bosonization, in Majid's
terminology) of $R$ and $H$.
 
\definition{Definition}Let $A$ be a Hopf algebra whose coradical $A_{0}$ is
a Hopf subalgebra. The braided Hopf algebra $R$ described above shall be
called the {\it diagram} of $A$. \enddefinition

\bigpagebreak
The general principle we propose is as
follows: first we analyze the diagram $R$ of $A$,
then we transfer the information to $\gr A$ by bosonization,
finally we lift it from $\gr A$ to $A$ via the filtration.

\bigpagebreak
$R$ is a {\it graded} braided Hopf  algebra and its
coradical is trivial: $R_{0} = R(0) = k1$. We denote by $P(R)$ the
space  of primitive elements of  $R$.  We see, considering the
coradical filtration, that   $P(R) \neq 0$,  because $\dim R >1$; this
last condition just means that $A$ is not cosemisimple.  In other
words,  the Hopf algebras $R$ we need to study are of a very
special kind.

 The first natural examples of such braided Hopf algebras
are the well-known  quantum linear spaces.  We give a
characterization of finite dimensional quantum linear spaces in
Section 3; see Proposition 3.5. 

If $\Gamma$ is a finite abelian group,
a quantum linear space over $\Gamma$ is given by elements $g_{1}, \dots, g_{\theta} \in \Gamma$, and characters $\chi_{1}, \dots,
\chi_{\theta} \in \widehat{\Gamma}$ satisfying

\roster \item" "  $q_{i}:= \chi_{i} (g_{i}) \neq 1$, \quad for all $i$,

\medpagebreak
\item" " $\chi_{j} (g_{i}) \chi_{i} (g_{j}) = 1$, \quad for all  $i
\ne j$. 
\endroster

The quantum linear space $\Cal R$ = $\Cal R(g_1, \dots,
g_{\theta}; \chi_1, \dots, \chi_{\theta})$ is then the braided Hopf algebra over $k\Gamma$ generated as an algebra by primitive elements  $x_{1},\dots, x_{\theta}$,
with relations 

\roster \item " "$x_{1}^{N_{1}} = 0, \dots,
x_{\theta}^{N_{\theta}} = 0,$ 

\medpagebreak 
\item " "$x_{i}x_{j} = \chi_{j}(g_{i}) x_{j}x_{i}$, if $i \ne j$.
\endroster

 The elements $x_i$ are $g_i$-graded and the action of $\Gamma$ on $x_i$ is via the character $\chi_i$. To each such quantum linear space we define a compatible datum $\Cal {D}$ consisting of scalars $\mu_{i} \in \{0, 1\}$ for each $i$, $1\leq i \leq
\theta$, and
 $\lambda_{ij} \in k$ for each $i, j$, $1\leq i < j \leq
\theta$ satisfying conditions (5.1) and (5.2). We define for each such datum
a pointed Hopf algebra $\Aq$ in section 5.   Then we prove our main result

\proclaim{Lifting Theorem 5.5} Let $\Cal R$ = $\Cal R(g_1, \dots, g_{\theta};
\chi_1, \dots, \chi_{\theta})$ be a quantum linear space over the finite
abelian group $\Gamma$. Then pointed Hopf algebras $A$ with coradical
$k(\Gamma)$ and diagram $\Cal R$ are exactly of the form $\Aq$ for some
compatible datum $\Cal D$. \endproclaim

\bigpagebreak

Let $p$ be an odd prime number and let $\Bbb G_p$ denote the group of
$p$-th roots of 1 in $k$.
 
As an application of the Lifting Theorem, we classify pointed Hopf algebras of dimension $p^3$.


 A Hopf algebra of dimension $p$ is isomorphic to a group algebra
by Zhu's theorem \cite{Z}.

The only pointed non-cosemisimple Hopf algebras of dimension
$p^2$ are the Taft  algebras $T_k(q) = T(q)$, $q \in \Bbb G_p - 1$;
see Section 1.

In dimension $p^3$, we have the following list of pointed
non-cosemisimple  Hopf algebras over $k$, for each $q\in \Bbb
G_p - 1$:

\bigpagebreak

 \roster  \item"(a)" The product Hopf algebra $T(q)\otimes k(\Bbb
Z/(p))$.

\bigpagebreak \item"(b)" The Hopf algebra $\widetilde{T(q)} :=
k\langle g, x\vert \; gxg^{-1} = q^{1/p}x, \; g^{p^2} = 1, \; x^p =
0\rangle$. Here     $q^{1/p}$ is a $p$-th root of $q$. Its
comultiplication is determined by $\Delta(x) = x\otimes g^p + 1
\otimes x$, $\Delta(g) = g\otimes g$. 

\bigpagebreak\item"(c)" The Hopf algebra $\widehat{T(q)} :=
k\langle g, x\vert  \; gxg^{-1} = qx,  \;g^{p^2} = 1,  \;x^p =
0\rangle$.
 Its comultiplication is determined by $\Delta(x) = x\otimes g + 1
\otimes x$, $\Delta(g) = g\otimes g$.

\bigpagebreak\item"(d)" The Hopf algebra $\bold r(q) := k\langle
g, x\vert  \; gxg^{-1} = qx,  \;g^{p^2} = 1, \; x^p = 1 - g^p\rangle$.
 Its comultiplication is determined by $\Delta(x) = x\otimes g + 1
\otimes x$, $\Delta(g) = g\otimes g$.

\bigpagebreak\item"(e)" The Frobenius-Lusztig kernel $\bold
u(q) :=   k\langle g, x, y\vert  \; gxg^{-1} = q^2x,  \; gyg^{-1} =
q^{-2}y, 
 \;g^{p} = 1,\; x^p = 0,  \; y^p = 0,  \; xy-yx = g-g^{-1}\rangle.$
 Its comultiplication is determined by $\Delta(x) = x\otimes g + 1
\otimes x$, $\Delta(y) = y\otimes 1 + g^{-1} \otimes y$,
$\Delta(g) = g\otimes g$.

\bigpagebreak\item"(f)" For each $m\in \Bbb Z/(p) -0$,  the book
Hopf algebra $\hq :=   k\langle g, x, y\vert  \; gxg^{-1} = qx,  \;
gyg^{-1} = q^{m}y, \; g^{p} = 1,  \; x^p = 0, \;  y^p = 0,  \; xy-yx =
0\rangle.$
 Its comultiplication is determined by $\Delta(x) = x\otimes g + 1
\otimes x$, $\Delta(y) = y\otimes 1 + g^{m} \otimes y$, $\Delta(g)
= g\otimes g$. \endroster

We prove \proclaim{Theorem 0.1} Any
non-cosemisimple  pointed Hopf algebra of order $p^3$ 
 is isomorphic to one  in the list above.\endproclaim

\bigpagebreak

 We shall see that there are no isomorphisms
between different   Hopf algebras in the list, except for book
algebras, where $\hq$ is isomorphic to $\bold h(q^{-m^{2}},
m^{-1})$; {\it cf.} Section 1.

\bigpagebreak
Let us define the {\it index} of a Hopf subalgebra $H$ of a finite
Hopf algebra $A$ as the ratio $\dim A/ \dim H$; it is an integer
because of the Theorem of Nichols-Zoeller  \cite{NZ}.

\medpagebreak
Theorem 0.1 is a consequence of  Theorem 0.2:

\proclaim{Theorem 0.2}Let $H = k(\Gamma)$, where $\Gamma$ is a
finite non-trivial abelian group; say $\Gamma =  \langle
y_{1}\rangle \oplus \dots \langle y_{\sigma}\rangle$, $y_{\ell} \neq 0$.
Let $M_{\ell}$ denote the order of $y_{\ell}$, $1\le \ell\le \sigma$.  Let $A$
be a pointed Hopf algebra with coradical  $H$.

\medpagebreak (A). Assume that the index of $H$ in $A$ is $p$.
Then there exist $g\in \Gamma$ and a character $\chi \in
\widehat{\Gamma}$ such that $q := \chi(g)$ has order $p$ and $A$
can be represented by generators $h_{\ell}$, $1\le \ell\le \sigma$, $a$, and
relations 

\roster \item"(0.2)" $h_{\ell}h_{t} = h_{t}h_{\ell}$, $h_{\ell}^{M_{\ell}} = 1$
for all $1\le \ell, t\le \sigma$
\item"(0.3)" $a^{p} = \mu(1- g^{p})$, with
$\mu$ either 0 or 1; \item"(0.4)" $h_{\ell}ah_{\ell}^{-1} = \chi(y_{\ell})a$, 
for all $1\le \ell\le \sigma$.
 \endroster The Hopf algebra structure of $A$ is
determined by  $$\Delta(h_{\ell}) = h_{\ell}\otimes h_{\ell}, \qquad \Delta(a) =
a\otimes 1 + g\otimes a,\quad   1\le \ell\le \sigma.$$

\medpagebreak  Assume that the index of $H$ in $A$ is $p^{2}$.
Then there are two possibilities:

(B$_{1}$). There exist $g\in \Gamma$ and a character $\chi \in
\widehat{\Gamma}$ such that $q := \chi(g)$ has order $p^{2}$ and
$A$ can be presented by generators $h_{\ell}$, $a$, $1\le \ell\le \sigma$, and
relations (0.2),

\roster  \item"(0.5)" $a^{p^{2}} = \mu(1- g^{p^{2}})$, with $\mu$
either 0 or 1; \item"(0.6)" $h_{\ell}ah_{\ell}^{-1} = \chi(y_{\ell})a$, $1\le \ell
\le \sigma$. \endroster The Hopf algebra structure of $A$ is determined by 
$$\Delta(h_{\ell}) = h_{\ell}\otimes h_{\ell}, \qquad \Delta(a) = a\otimes 1 +
g\otimes a,\quad  1\le \ell\le \sigma$$.

\medpagebreak (B$_{2}$). There exist $g_{1}, g_{2}\in \Gamma$
and  characters $\chi_{1}, \chi_{2} \in \widehat{\Gamma}$ such that
$q_{1} := \chi_{1}(g_{1})$ and $q_{2} = \chi_{2}(g_{2})$ have order
$p$, $\chi_{1} (g_{2}) \chi_{2} (g_{1}) = 1$ and $A$ can be
presented by generators $h_{\ell}$, $a_{i}$, $1\le \ell  \le \sigma$, $i =
1,2$, and relations (0.2),

\roster  \item"(0.7)" $a_{i}^{p} = \mu_{i}(1- g_{i}^{p})$, with
$\mu_{i}$ either 0 or 1, $i = 1,2$;  \item"(0.8)"  $h_{\ell}a_{i}h_{\ell}^{-1} =
\chi_{i}(y_{\ell})a_{i}$, $1\le \ell \le \sigma$, $i =
1,2$; \item"(0.9)" $a_{1}a_{2} - \chi_{2}(a_{1}) a_{2}a_{1} =
\lambda(1 - g_{1}g_{2})$,  with $\lambda$ either 0 or 1. \endroster
If $\lambda \neq 0$, then $\chi_{1}\chi_{2} = 1$. The Hopf algebra
structure of $A$ is determined by  
$$\Delta(h_{\ell}) = h_{\ell}\otimes h_{\ell},
\qquad \Delta(a_{i}) = a_{i}\otimes 1 + g_{i}\otimes a_{i},\quad i =
1,2, \quad 1\le \ell\le \sigma.$$

\endproclaim
The proof of Theorem 0.2 follows from the above principle: we
show, via the mentioned characterization, that a braided Hopf
algebra of our special type and of dimension $p$ or $p^{2}$ is
necessarily  a quantum linear space (Lemma 5.6). We then deduce
Theorem 0.2 from the Lifting Theorem (5.5).

\bigpagebreak
We shall give more  applications of this principle in a separate
article. We shall generalize  the  basic Theorem of Taft and
Wilson (\cite{TW}, \cite{M, Thm. 5.4.1}) to the case of Hopf algebras
whose coradical is a Hopf subalgebra. This Theorem is the key
point in the proof of the following result (see {\it e. g.} \cite{N, p.
1545}, \cite{AS2, Prop. 3.1}): If $A$ is a pointed non-cosemisimple
finite dimensional Hopf algebra, with coradical  $k(\Gamma)$
where $\Gamma$ is abelian, then there exist  $g\in \Gamma$, a 
$k$-character $\chi$ of $\Gamma$ such that $\chi(g) \ne 1$, and
$x\in A$, $x \notin k(\Gamma)$ such that  $$hxh^{-1} = \chi(h)x
\quad \forall h \in \Gamma,\quad \Delta(x) = x\otimes g +
1\otimes x. $$  The preceding statement is the initial
point in existing attempts of classifications of various kinds of
pointed Hopf algebras.

\bigpagebreak 
The Lifting Theorem 5.5 has an extra bonus. In 1975, Kaplansky
formulated a series of conjectures on Hopf algebras. Under the hypothesis
that the characteristic of the ground field does not divide the positive
integer $n$, one of these conjectures states that there are only a finite
number (up to isomorphism) of Hopf algebras of  dimension $n$. In this
direction the  following result was proved by Stefan: 
The set of types of semisimple and cosemisimple Hopf algebras of a
given dimension is finite (in any characteristic). See \cite{St}; a more direct
proof (showing at the same time finiteness of the number of automorphism and right coideal subalgebras) is given in \cite{S}.  

Our Lifting Theorem easily produces counterexamples to Kaplansky's conjecture.

\proclaim{Theorem 0.3} There exist an infinite family of non-isomorphic
pointed Hopf algebras of order $p^{4}$. \endproclaim

\bigpagebreak Let us say that a Hopf algebra is {\it very simple} if
\roster\item"(i)"  it has no non-trivial normal Hopf subalgebra, and
 \item"(ii)"   it can not be constructed by bosonization in a
non-trivial  way. \endroster

Usually, a Hopf algebra is called simple if it satisfies only (i); see
for instance \cite{A}. However, Taft algebras and book algebras
are simple in this sense- this follows from the criteria in
\cite{AS1}; but they are  analogues of solvable algebraic groups
and it is hard to accept their  simplicity.  On the other hand,
bosonization is also a mean to build Hopf algebras from  smaller
ones- though one of them is a braided Hopf algebra.  Also, Taft
and book algebras can be build by bosonization. For these
reasons, we propose this new definition.

\bigpagebreak  Theorem 0. 1 has the following consequence:

\proclaim{Corollary 0.4} The only   pointed Hopf algebras of order
$p^3$  which are very simple, are the Frobenius-Lusztig kernels 
$\bold u(q)$ of type $A_{1}$.\endproclaim

The Corollary follows from the considerations in Section 1. So far,
the only known very simple Hopf algebras of order $p^3$ are the 
Frobenius-Lusztig kernels $\bold u(q)$ and their duals; see {\it
1.7}.

Let us briefly indicate the contents of the paper. In Section 1, we
give some information about the Hopf algebras mentioned above.
Section 2 is devoted to basic facts supporting the principle. 
 In Section 3 we discuss finite dimensional quantum linear spaces. In
Section 4, we discuss possible quantum linear spaces over abelian groups. 
Theorem 0.2,  respectively Theorem 0.3,  Theorem 0.1,  are proved in Section
5,  respectively Section 6, Section 7.

Theorem 0.1 of this paper was announced at the Colloque "K-theory, cyclic
homology and group representations", CIRM, Luminy,  (july 1997); and at the
"XLVI Reuni\'on de Comunicaciones Cient\i ficas de la Uni\'on Matem\'atica 
Argentina", C\'ordoba, (september 1997), where also a counterexample to
Kapklansky's conjecture was described. The list appears already in the
preprint version  of \cite{AS2}, Trabajos de Matem\'atica 42/96, FaMAF.

Just before submitting this paper we learned that Theorem 0.1 resp. Theorem 0.3
also appears in recent unpublished work of S. Caenenepeel, S. Dascalescu resp.  M. Beattie, S. Dascalescu L. Gr\"unenfelder, and also S. Gelaki. The methods of these authors seem to be quite different from ours.\footnote{This paper was submitted on November 1997.}

\subsubhead Conventions\endsubsubhead
If $C$ is a coalgebra, we denote  by $G(C)$ the set of group-like
elements of $C$. If $g$, $h \in G(C)$, then we denote $P_{g,h}(C) =
P_{g,h} = \{x\in C: \Delta(x) = x\otimes h + g \otimes x\}$; the
elements of $P_{g,h}$ are called skew-primitives.  When $B$ is a
bialgebra, $P_{1,1}(B)$ is just the space $P(B)$ of primitive
elements.

If $A$ is an algebra, $\Alg(A, k)$ denotes the set
of all algebra maps from $A$ to $k$.

If $\Gamma$ is a group, we denote by $\widehat{\Gamma}$ the
group of characters (one-dimensional representations over $k$)
of $\Gamma$.

\subhead \S 1. About the Hopf algebras in the list\endsubhead 

\subsubhead  1.0\endsubsubhead Let $\xi\in k$ be a root of 1 of 
order $N \ge 2$.
 The  Taft algebra $T_k(\xi) = T(\xi)$ is the algebra   $ k\langle g,
x\vert \; gxg^{-1} =\xi x, \; g^{N} = 1, \; x^N = 0\rangle$. Its  Hopf
algebra structure is given by $$\aligned \Delta(g) = g\otimes g,
\quad \Cal S(g) &= g^{-1}, \quad \varepsilon(g) = 1, \\ \Delta(x) =
x\otimes g + 1\otimes x, \quad \Cal S(x) &= -xg^{-1},  \quad
\varepsilon(x) = 0. \endaligned$$ The dimension of $T(\xi)$ is
$N^{2}$. It is known that $T(\xi)\simeq T(\xi)^*$ as Hopf algebras,
and that  $T(\xi)\simeq T(\widetilde \xi)$ only if $\xi = \widetilde
\xi$. 

A proper Hopf subalgebra $A$ of $T(\xi)$ is contained in $k[g]$:
this follows  easily looking at the coradical filtration of $A$. 
Therefore, if $A$ is a proper  Hopf subalgebra or quotient of
$T(\xi)$, then the order of $A$ divides $N$.

Semisimple Hopf algebras of order $p^2$ are group algebras
\cite{Ma}. The only pointed non-cosemisimple Hopf algebras of
order $p^2$ are the Taft algebras; a more precise
characterization of Taft algebras is given in  \cite{AS2}. In fact,
Taft algebras and group algebras are the only known Hopf
algebras of order $p^2$.

\bigpagebreak
\subsubhead  1.1\endsubsubhead The pointedness of the Hopf
algebras in the list is a consequence of the criteria \cite{M,
Lemma 5.5.1}. As in the proof of \cite{M, Lemma 5.5.5}, we
conclude that a Hopf algebra in the list has  coradical
$k(\Gamma)$, where $\Gamma$ is:
\roster \item $ \Bbb Z/(p) \times \Bbb Z/(p)$, in case (a);
\item $ \Bbb Z/(p^{2})$, in cases (b), (c), (d);
\item $ \Bbb Z/(p)$, in cases (e), (f).
\endroster
In particular, this is a first step towards deciding the
non-existence of isomorphisms between the different cases.

\bigpagebreak It is not difficult to see that all the Hopf algebras
in the list have dimension $p^{3}$; {\it e.g.}, using the Diamond
Lemma. Alternatively, say in case (b), let $A$ be a vector space
with a basis $g^{i}x^{j}$, $0\le i \le p^{2}-1$, $0\le j \le p-1$. It is
possible to write down explicitly a multiplication table for $A$
such that the defining relations hold; $A$ is then an associative
algebra. Hence there is an epimorphism $\widetilde{T_{k}(q)} \to
A$. But it is easy to see that $\widetilde{T_{k}(q)}$ has dimension
at least $p^{3}$; therefore the dimension is $p^{3}$. This idea 
applies to the other cases as well.

\bigpagebreak \subsubhead  1.2\endsubsubhead The Hopf algebra
$\widetilde{T_{k}(q)}$  does not depend, modulo isomorphisms,
 upon the choice of the  $p$-th root of $q$.  Indeed, let
$\widetilde{T_{j}} :=  k\langle h, y\vert \; hyh^{-1} = q^{1/p + j}y, \;
h^{p^2} = 1, \; y^p = 0\rangle$, with comultiplication $\Delta(h) =
h\otimes h$, $ \Delta(y) = y\otimes h^{p}+ 1\otimes y.$ Then one
has an isomorphism of Hopf algebras $\widetilde{T_{k}(q)} \to
\widetilde{T_{j}}$ determined by $x \mapsto y$, $g \mapsto h^{1 -
pj}$.

Notice that $\widetilde{T_{k}(q)}$ is a cocentral extension of $
k(\Bbb Z/p)$ by  a Taft algebra: $$1 @>>>T_{k}(q) @>>>
\widetilde{T_{k}(q)} @>>> k(\Bbb Z/p)  @>>> 1.$$

\bigpagebreak

\subsubhead  1.3\endsubsubhead  The Hopf algebra
$\widehat{T_{k}}(q)$  is dual to $\widetilde{T_{k}(q)}$. It is a
central extension of a Taft algebra: $$1 @>>> k(\Bbb Z/p) @>>>
\widehat{T_{k}(q)} @>>> T_{k}(q) @>>> 1;$$ where the central Hopf
algebra is generated by $g^p$. 
 It is clear that no group-like element of $\widetilde{T_{k}(q)}$ is
central;  hence $\widetilde{T_{k}(q)}$ and  $\widehat{T_{k}(q')}$
can not be isomorphic for any $q, q'$.

\bigpagebreak \subsubhead  1.4\endsubsubhead The Hopf algebra
$\bold r(q)$ is also a central extension of a Taft algebra: $$1 @>>>
k(\Bbb Z/p) @>>> \bold r(q) @>>> T_{k}(q) @>>> 1;$$ again, the
central Hopf algebra is generated by $g^p$. This Hopf algebra was first
considered by Radford \cite{R1}. The dual Hopf algebra  $\bold r(q)^*$ is not
pointed-- see {\it loc. cit.};  hence cases (b), (c) and (d) have no
intersection. 

\bigpagebreak In all three cases, the Hopf algebras are not
isomorphic for different values of $q$. This can be shown via the
first term of the coradical filtration. Indeed, it is enough to
consider $\widehat{T_{k}(q)}$, since $\widetilde{T_{k}(q)}\simeq
\left(\widehat{T_{k}(q)}\right)^{*}$ and $\gr \bold r(q) \simeq
\widehat{T_{k}(q)}$.

\bigpagebreak \subsubhead  1.5\endsubsubhead The
Frobenius-Lusztig kernel $\bold u(q)$ is the simplest example of
the finite dimensional Hopf algebras introduced in  \cite{L1},
\cite{L2}. It is easy to see that it has no non-trivial
representation of dimension 1. Looking at its coradical filtration
\cite{T}, we conclude that  $\bold u(q)$ and $\bold u(q')$ are not
isomorphic unless $q = q'$. It is not difficult to see that $\bold
u(q)$ has no non-trivial quotient Hopf algebra \cite{T}; hence it is
very simple. See also \cite{AS1}.

\bigpagebreak \subsubhead  1.6\endsubsubhead Information
about  book Hopf algebras  can be found in \cite{AS2, \S 6};  $\bold
h(q,p-1)$ was already considered in \cite{R2, p. 352}- without
assuming that the order of $q$ is prime. $\hq$ and $\bold
h(\widetilde q, \widetilde m)$ are  isomorphic if, and only if,
$(\widetilde q, \widetilde m) = (q,m)$ or $(q^{-m^{2}}, m^{-1})$ 
\cite{AS2, Prop. 6.5}. The dual Hopf algebra $(\hq)^{*}$ is
isomorphic to $\bold h(q, -m)$ \cite{AS2, Prop. 6.7}; in particular, 
$\hq$ has $p$ different representations of dimension 1 and hence
types (e) and (f) have no intersection.

Book algebras can be obtained by bosonization \cite{AS2}; see
also Section 3. By the criteria in \cite{AS1}, a book algebra is
simple.

\bigpagebreak

\subsubhead  1.7\endsubsubhead Semisimple Hopf algebras of
order $p^3$ were classified  by Masuoka \cite{Ma}: there are $p+8$
isomorphism types, namely three group algebras of abelian
groups;  two group algebras of non-abelian groups and their duals;
$p + 1$ non-commutative, non-cocommutative  Hopf algebras
constructed by extension. 

In addition to the already mentioned Hopf algebras of order $p^3$
there are also  the dual Hopf algebras   $(\bold u(q))^*$ and $(\Cal
R(q))^*$.  Among all these  Hopf algebras of order $p^3$, only the
Frobenius-Lusztig kernels and their duals are very simple in the
sense of the Introduction.

\subhead \S 2. The coradical filtration and the  associated graded
Hopf algebra\endsubhead 

\subsubhead  2.0\endsubsubhead Let $B$, $H$ and $R$ be as in (0.1).
Then $R$ is a braided Hopf algebra in the category 
${}_H^H\Cal{YD}$ of Yetter-Drinfeld modules over $H$. See
\cite{R3}, \cite{Mj}. We recall the explicit form of this structure;
we follow the conventions of \cite{AS2}. 

\medpagebreak
The action of $H$  on $R$ is given by the adjoint representation
composed with $\gamma$. The coaction is $(\pi\otimes
\id)\Delta$.  These two structures are related by the
Yetter-Drinfeld condition: $$\delta_{R}(h.r) = h_{(1)} r_{(-1)}
\Ss(h_{(3)}) \otimes h_{(2)} .r_{(0)}.   $$ Hence, $R$ is an object  of
${}_H^H\Cal{YD}$.

\medpagebreak
Moreover, $R$ is a subalgebra of $B$ and a coalgebra with
comultiplication  $$\Delta_{R}(r) = r_{(1)} \gamma\pi\Ss(r_{(2)})
\otimes r_{(3)};$$ the counit is the restriction of the counit of $B$. 
To avoid confusions, we denote here the comultiplication of $R$ in
the following way: $$\Delta_{R}(r) = \sum r^{(1)}\otimes r^{(2)},$$ or
even we omit sometimes the summation sign.

\medpagebreak
The multiplication $m$ and the comultiplication $\Delta$ of $R$
satisfy $$\Delta m =  (m \otimes m)(\id \otimes c \otimes \id)
(\Delta\otimes \Delta).$$ Here $c$ is the commutativity constraint 
of ${}_H^H\Cal{YD}$; explicitly $$c_{M, N}(m\otimes n) = m_{(-1)}. n
\otimes m_{(0)},$$ for $M, N \in {}_H^H\Cal{YD}$, $m\in M$, $n\in N$.

\medpagebreak The map $\Ss_{R}: R \to R$ given by $$\Ss_{R}(r) =
\gamma \pi(r_{(1)}) \Ss_{B} r_{(2)}$$ (where $\Ss_{B}$ is the
antipode of $B$) is the antipode of $R$, {\it i. e.} the inverse of the
identity in $\End R$ for the convolution product. Hence $R$ is a
braided Hopf algebra in  ${}_H^H\Cal{YD}$.

\medpagebreak
Conversely, given $H$ and a Hopf algebra $R$ in  ${}_H^H\Cal{YD}$,
the tensor product $B = R\otimes H$ bears a Hopf algebra
structure, denoted $R\# H$, via the smash product and smash
coproduct:
$$(r\#   h)(s\#   f) = r(h_{(1)}.s) \#   h_{(2)}f,
\qquad \Delta(r\#   h) = r^{(1)}\#  (r^{(2)})_{(-1)}h_{(1)} \otimes
(r^{(2)})_{(0)} \# h_{(2)}.  \tag 2.1$$

Let $\pi: R \# H\to H$ and $\gamma: H \to R \# H$ be the 
maps $$\pi(r\# h) = \varep(r)h, \quad \gamma(h) = 1\# h.$$ Then
$\gamma$ is a section of $\pi$ and we are in the setting (0.1). We
term $B = R\# H$, following Majid, the bosonization of $R$.

\subsubhead  2.1\endsubsubhead We recall that a graded Hopf
algebra is a Hopf algebra $G$ together with a grading $G =
\oplus_{n\ge 0} G(n)$ which is simultaneously an algebra and a
coalgebra grading \cite{Sw, Section 11.2}.  In particular,
$\varep(G(n)) = 0$ for  $n >0$ and the antipode is a homogeneous
map of degree 0.

It turns out that $G(0)$ is  a Hopf subalgebra of $G$ and that the
inclusion $\gamma: G(0) \hookrightarrow G$ is a section of the
projection $\pi: G \to G(0)$ with kernel $\oplus_{n\ge 1} G(n)$.  Let
$$R = \{a \in G: (\id \otimes \pi) \Delta(a) = a\otimes 1\}.$$ We
know that  $R$ is a braided Hopf algebra in
${}_{G(0)}^{G(0)}\Cal{YD}$ and  $G$ is the bosonization of $R$.

We shall say that  a braided  Hopf algebra with a grading of
Yetter-Drinfeld modules is {\it graded} if the grading  is
simultaneously an algebra and  coalgebra and  grading.

\proclaim{Lemma 2.1}Keep the notation above. 
\roster \item"(i)" $R$ is a graded subalgebra of $G$: $R =
\oplus_{n\ge 0} R(n)$, where $R(n) = R \cap G(n)$. 

\item"(ii)" With respect to this grading, it is a braided
graded Hopf algebra.

\item"(iii)" $G(n) = R(n)\# G(0)$ and $R_{0} = R(0) = k1$.
\endroster      \endproclaim

\demo{Proof}Let $r \in R$ and let us decompose $r = \sum_{j}
r_{j}$, with $r_{j} \in G(j)$. Then we write
$$\Delta_{G}(r_{j}) = \sum_{h} r_{j, h}\otimes r_{j}^{h},$$
where $r_{j, h}\in G(h)$, $r_{j}^{h}\in G(j-h)$.   Clearly, $\pi(r_{j}^{h})
= 0$ unless $h=0$. Hence  $(\id \otimes \pi) \Delta(r_{j}) \in
G(j)\otimes G(0)$. By definition of $R$,
$$\sum_{j} r_{j}\otimes 1 = \sum_{j} (\id \otimes \pi)
\Delta(r_{j});$$ taking homogeneous components, we see that
$r_{j}\otimes 1 =  (\id \otimes \pi) \Delta(r_{j})$, {\it i.e.} that
$r_{j}\in R$. This proves (i).

It follows from the definition of the action and coaction that each
$R(n)$ is a submodule and subcomodule. That is, $R = \oplus_{n\ge
0} R(n)$ is a grading in ${}_{G(0)}^{G(0)}\Cal{YD}$.

It is not difficult that $R$ is a graded coalgebra. Indeed, if $r\in
R(j)$ then we write
$$(\Delta_{G}\otimes \id)\Delta_{G}(r) = \sum_{h, t} a_{h}\otimes
b_{t}\otimes c_{j-t-h},$$
where $a_{h}\in G(h)$, $b_{t}\in G(t)$, $c_{j-t-h}\in G(j-t-h)$. Hence
$$\Delta_{R}(r) = \sum_{h, t} a_{h}\pi\Ss (b_{t}) \otimes
c_{j-t-h} = \sum_{h} a_{h}\pi\Ss (b_{0}) \otimes
c_{j-h} \in \oplus_{h} R(h) \otimes R(j-h).  $$
Now we prove (iii). The first claim is evident, since $G(n)
\supseteq R(n)\# G(0)$ and $G = \oplus_{n\ge 0}G(n) =
\oplus_{n\ge 0} R(n)\otimes G(0)$. As for the second, we know
that  $R_{0} \subseteq  R(0)$ and $ R(0) = k1$. Indeed, the
contention follows since the coradical is contained in the zero
part of any coalgebra filtration \cite{M, 5.3.4}; the equality
follows by definition. These two facts imply that $R_{0} = R(0) =
k1$.  \qed\enddemo

\subsubhead  2.2\endsubsubhead  Let $A$ be a Hopf algebra  and
assume that its coradical $A_0$ is a Hopf subalgebra (for
instance, $A$ is  pointed). Then the coradical filtration is in fact a
Hopf algebra  filtration and the associated graded algebra  $$  \gr
A = \oplus_{n \ge 0} \gr A(n) = \oplus_{n \ge 0} A_n/A_{n-1}$$
(with $A_{-1} = 0$) is a graded Hopf algebra.  See \cite{M, 5.2.8}. If
$A$ has finite dimension $N$,  then $\gr A$ also has dimension $N$.

\proclaim{Lemma 2.2}If $\gr A$ is generated as an algebra by $\gr
A(0) \oplus \gr A(1)$ then  $A$ is generated as an algebra by
$A_{0} + A_1$.\endproclaim 
\demo{Proof} This can be checked directly,  or via the following
argument: $A_1$, $\gr A(1)$ are $A_0$-bimodules,  and the
projection $A \to \gr A(1)$ is a bimodule homomorphism.  Since
$A_0$ is semisimple by \cite{LR}, $A_1 \simeq A_0 \oplus \gr
A(1)$ as $A_0$-bimodules. We can  consider the tensor algebra
$T_{A_0}(\gr A_1)$ and the corresponding map $\pi: T_{A_0}(\gr
A_1) \to A$- see \cite{N, Prop. 1.4.1}.  This map is compatible with
filtrations and the corresponding graded map is bijective. Then
$\pi$ is bijective \cite{B, \S 2, no. 8}. \qed\enddemo

\subsubhead  2.3\endsubsubhead Let $A$ be a Hopf algebra  whose
coradical $A_0$ is a Hopf subalgebra.

\proclaim{Lemma 2.3}The coradical filtration of $\gr A$ is given
by $$(\gr A)_m =  \oplus_{n \le m} \gr A(n). \tag 2.2$$ 
\endproclaim

\definition{Definition \cite{CM}} A graded coalgebra satisfying (2.2) is
called {\it coradically graded}.\enddefinition

\demo{Proof}We check this for $m= 0, 1$; the general case is
similar or else can be deduced  from these two cases by \cite{CM,
2.2}. 

First, $A_0 = \gr A(0) \subseteq (\gr A)_0$ because it is
cosemisimple.  Conversely, the filtration $\gr A(0) \subset \gr
A(0)\oplus \gr A(1) \subset \dots \subset  \oplus_{n \le m} \gr
A(n) \subset \dots $ is a coalgebra filtration hence $\gr A(0)
\supseteq (\gr A)_0$  \cite{Sw, 11.1.1}.

Now we consider $m=1$. Again,  $A_0\oplus A_1/A_0 \subseteq
(\gr A)_1$ is easy. Let $y\in \gr A$ and write $$y = \overline{y_0}
+  \overline{y_1} + \dots + \overline{y_m}, \quad  \overline{y_j}
\in A_j/A_{j-1}, \quad \overline{y_m}\neq 0.$$ Hence $$ \Delta(y)
= \sum_{j=0}^m \Delta(\overline{y_j}) \in \overline{y_m} +
\oplus_{r + s < m} \gr A(r) \otimes \gr A(s), $$ and $\Delta(
\overline{y_m}) = z_1 + z_2 + z_3$, with $ z_1  \in\gr A(m)
\otimes \gr A(0)$, $z_2 \in\gr A(0) \otimes \gr A(m)$,  $z_3 \in
\oplus_{r + s = m, r,s > 0} \gr A(r) \otimes \gr A(s).$

Now assume that $m>1$.  If $z_3 = 0$, $\overline{y_m} = 0$; and if 
$z_3 \neq 0$, $y\notin (\gr A)_1$. So $m$ should be 1 and
$A_0\oplus A_1/A_0  \supseteq (\gr A)_1$. \qed\enddemo

\subsubhead  2.4\endsubsubhead Let $G = \oplus_{n\ge 0} G(n)$ be
a coradically graded Hopf algebra. Let  $R$ be the associated
braided graded Hopf algebra, see {\it 2.1}.

\proclaim{Lemma 2.4}\roster \item"(i)" $R_{0} = k1 = R(0)$ and
$P(R) = R(1)$. \item"(ii)" $R$ is a  coradically graded coalgebra.  
\item"(iii)"  $G_{1} = G(0) \oplus
\left[ P(R)\# G(0)\right]$. \endroster
\endproclaim

\demo{Proof} (i). We know that  $R_{0} = R(0) = k1$ from Lemma
2.1. 

Let $r \in R(1)$. Then $\Delta_{R}(r) = r_{1}\otimes 1 + 1\otimes
r_{2}$, for some $r_{1}, r_{2} \in R(1)$. Applying $\id\otimes
\varep$ and $\varep\otimes \id$ to both sides of this equality, we
conclude that $r_{1} = r_{2} = r$. That is, $P(R) \supseteq R(1)$.

Let now $r \in P(R)$. If $\delta (r) = r_{(-1)} \otimes r_{(0)} \in
G(0)\otimes R$, then
$$\Delta_{G}(r)  = r\otimes 1 + r_{(-1)} \otimes r_{(0)}.$$ As
$G_{0} = G(0)$, we deduce that $r\in G_{1}$. But by hypothesis,
$G_{1} = G(0)\oplus G(1)$. We can assume that $r$ is
homogeneous; since $r\in R(0)$ is not possible, we see that $r \in
G(1)$. That is, $P(R) \subseteq R(1)$.

(ii). By \cite{CM, 2.2}, it is enough to consider the cases $m = 0, 1$.
The case $m=0$ is covered by (i). For $m=1$, we have, again by
(i), $$R_{1} = k1\oplus P(R) = R(0) \oplus R(1).$$

(iii).  This follows from Lemma 2.1 and (ii). \qed\enddemo

\subsubhead  2.5\endsubsubhead Let $H$ be a cosemisimple Hopf
algebra. Let  $R$ be a braided  graded Hopf algebra in the
category  ${}_H^H\Cal{YD}$. Let $G = R\# H$; it is easy to see that
$G$ is a graded Hopf algebra.

\proclaim{Lemma 2.5}If $R_{0} = k1 = R(0)$ and $P(R) = R(1)$, then
$G$ is a coradically graded Hopf algebra. \endproclaim

\subhead \S 3. Quantum linear spaces\endsubhead 

 As mentioned in the Introduction,  we are interested in  braided
Hopf algebras $R$ in the category  ${}_H^H\Cal{YD}$ of
Yetter-Drinfeld modules  over a Hopf algebra $H$. We use in this
section the notation of  \cite{AS2, Section 4}.

 A version of the following result appears in \cite{N, p. 1538}.

\proclaim{Lemma 3.1} Let $H$ be a finite dimensional Hopf
algebra. 

(i). Let $x \in H$ such that $\Delta(x) = x\otimes 1 + g\otimes x,$ 
$gx = xg$, for some $g \in G(H)$. Then $x$ is a scalar multiple of
$g-1$.

(ii). Let $R$ be a finite dimensional braided Hopf algebra in 
${}_H^H\Cal{YD}$. Let $x \in P(R)$ be a non-zero primitive element
such that  $\delta (x) = g\otimes x$,  $h. x = \chi(h)x$,  for some
$g \in G(H)$, $\chi \in \Alg(H, k)$ and for all $h\in H$.  Then $q:=
\chi(g) \neq 1$.\endproclaim 

\demo{Proof} Let $S$ be the subalgebra  of  $H$  generated by $g$
and $x$; by hypothesis, it is a commutative Hopf subalgebra and
hence it is cosemisimple by the Cartier-Kostant theorem. 
Looking at the expression of $x$ in terms of the decomposition of
$S$ in simple subcoalgebras, one concludes that $x = \lambda
(g-1)$, for some   $\lambda\in k$. This shows (i).

For (ii), we apply (i) to the element $x\# 1$ of $A = R\# H$; by (2.1), 
$\Delta(x) = x\otimes 1 + g\otimes x$ and $gx = qxg$. If $q=1$ then
$x = \lambda (g-1)$, for some   $\lambda\in k$. This implies $x = 0$, a
contradiction.  \qed\enddemo

\bigpagebreak Let $K$ be an arbitrary Hopf algebra. Let $g\in
G(K)$,   $\chi\in \Alg(K,k)$  such that $\chi(h) g =  h_{(1)}
\chi(h_{(2)}) g \Ss(h_{(3)})$, for all $h\in K$.  Let $N$ be the order
of $q :=\chi(g)$; we assume $N$ is finite.

Let $R = k[y] / (y^{N})$. Then $R$ is a  braided Hopf algebra in 
${}_K^K\Cal{YD}$ with $K$-module    and $K$-comodule structures
given by $$h.y^{t} = \chi^{t}(h)y^{t}, \qquad \delta_{R}(y^{t}) =
g^{t}\otimes y^{t},$$ and  comultiplication uniquely determined by
$\Delta_{R}(y) = y\otimes 1 + 1\otimes y$. This braided Hopf
algebra will be denoted $\Cal R(g, \chi)$.  The braided Hopf
algebras $\Cal R(g, \chi)$ and $\Cal R(\widetilde g, \widetilde
\chi)$ are isomorphic  only if $g = \widetilde g$ and $\chi =
\widetilde \chi$. See \cite{AS2, Lemma 8.1}.

\proclaim{Theorem       3.2} Let $H$ be a finite dimensional
semisimple Hopf algebra. Let $R$ be a finite dimensional braided
Hopf algebra in  ${}_H^H\Cal{YD}$. Assume that \roster  \item 
$R_{0} = k1$, where $R_0$ is the coradical of $R$. \item  $ \dim
P(R)  = 1$.
 \endroster Then there exist  $g\in G(H)$, and  $\chi\in \Alg(H, k)$
such that $R \simeq \Cal R(g, \chi)$. \endproclaim \demo{Proof} By
(1), $P(R) \neq 0$. As $P(R)$ is a Yetter-Drinfeld submodule of $R$,
condition (2) implies the existence of $g\in G(H)$ and a character
$\chi\in H^{*}$ such that  $$\delta (x) = g\otimes x, \qquad h. x =
\chi(h)x \quad \forall h\in G(H), x \in P(R).$$   Let $N$ be the order
of $q = \chi(g)$. Fix $x \in P(R)$, $x\ne 0$. By the quantum binomial
formula, $\dim k[x] = N$ and $x^{N} = 0$. In fact, $ k[x]  \simeq \Cal
R(g, \chi)$. Hence we only need to prove that $R =  k[x]$.

Consider the algebra $R^{*}$.  It has a unique maximal ideal,
namely $\Cal M := R_{0}^{\bot}$. $\Cal M$, and {\it a fortiori\/} $\Cal
M^{2}$ and $\Cal M / \Cal M^{2}$, are Yetter-Drinfeld modules.
Observe that $\Cal M / \Cal M^{2} \simeq \left(P(R) \right)^{*}$ as
$H$-modules. But, since $H$ is semisimple, the projection $\Cal M
\to \Cal M / \Cal M^{2}$ has an $H$-linear section. Whence there
exists $T \in \Cal M - \Cal M^{2}$ such that  $$\delta (T) =
g^{-1}\otimes T, \qquad h. T = \chi^{-1}(h)T \quad \forall h\in H.$$
It is not difficult to show that  $R^{*} = k[T]$. Hence $1, T,  T^{2},
\dots, T^{d-1}$ is a basis of $R^{*}$, where $d = \dim R$, and we
can  consider its dual basis $t_{0}, t_{1}, \dots, t_{N}, \dots$ in $R$.
Note that  $\delta (T^{j}) = g^{-j}\otimes T^{j}$,  $h. T^{j} =
\chi^{-j}(h)T^{j}$, $\forall h\in H$;  hence $\delta (t_{j}) =
g^{j}\otimes t_{j}$, $ h. t_{j} = \chi^{j}(h)t_{j}$, $\forall h\in H$.

On the other hand, consider the coradical filtration of $R$:
$$R_{0}= k1 \; \subseteq \; R_{1} = k1 \oplus P(R) \;\subseteq \;
\dots \;\subseteq\; R_{j} \;\subseteq\; \dots$$ If $j \le N - 1$, $1,
x, \dots, x^{j}$ belong to $R_{j}$. As  $(R_{j})^{*}  \simeq R^{*}/\Cal
M^{j+1}$, we conclude that $1, x, \dots, x^{j}$ form a basis of 
$R_{j}$. Hence there exist $\lambda_{j} \in k$ such that
$$\lambda_{j}t_{j} = x^{j}, \qquad j \le N - 1.$$ Now assume $d > N$
and let $z = t_{N}$. Then $$\Delta(z) = \Delta(t_{N}) = \sum_{0 \le i
\le N} t_{i}\otimes t_{N - i} = z\otimes 1 + 1 \otimes z + \sum_{1
\le i \le N -1} \lambda_{i}\lambda_{N - i }x^{i}\otimes x^{N - i}.$$ 
Therefore the subalgebra $k\langle x, z\rangle$ is a Hopf subalgebra of $R$. Now
let us compute  $$\multline \Delta(xz) = (x\otimes 1 + 1 \otimes
x) \left( z\otimes 1 + 1 \otimes z + \sum_{1 \le i \le N -1}
\lambda_{i}\lambda_{N - i }x^{i}\otimes x^{N - i}\right) \\ =
xz\otimes 1 + x \otimes z + \sum_{1 \le i \le N -1}
\lambda_{i}\lambda_{N - i }x^{i+1}\otimes x^{N - i} \\ +  z\otimes
x  + 1 \otimes zx + \sum_{1 \le i \le N -1} \lambda_{i}\lambda_{N -
i }q^{i}x^{i}\otimes x^{N + 1 - i}. \endmultline$$ Hence $xz\in
R_{N+1}$; similarly, also $zx\in R_{N+1}$.  We conclude, looking  at
the decomposition in $H$-submodules, that  $$xz = a t_{N+1} + bx,
\qquad zx = c t_{N+1} + dx,$$ for some $a,b,c,d\in k$. It follows
from this that $k[x]$ is a normal Hopf subalgebra of $k\langle x, z\rangle$, and
we can form the quotient Hopf algebra. The image of $z$ in this
finite dimensional braided Hopf algebra is invariant and primitive,
hence 0.  Then $k[x] = k\langle x, z\rangle$, a contradiction. \qed\enddemo

\bigpagebreak Let $K$ be an arbitrary Hopf algebra.

\definition{Definition}We shall say that a braided Hopf algebra $R$
in ${}_K^K\Cal{YD}$ satisfies {\it hypothesis (A)} if there exist a
basis $x_{1}, \dots, x_{\theta}$ of $P(R)$, and $g_{1}, \dots,
g_{\theta} \in G(K)$, $\chi_{1}, \dots, \chi_{\theta} \in \Alg(K,k)$
such that  for all $j$, $1\le j \le \theta$,

\roster \item" " $\delta (x_{j}) = g_{j}\otimes x_{j}$, \quad
 $h. x_{j} = \chi_{j}(h)x_{j}$, \quad for all $h\in K$,  and

\smallpagebreak
 \item" " the order $N_j$ of $q_{j}:= \chi_{j} (g_{j})$ is finite.
\endroster\enddefinition

If $K = k(\Gamma)$ is the group algebra of a  finite abelian group
$\Gamma$, hypothesis (A) always holds.  

If $q \in k$ and $0 \le i \le n < \ord q$, we set $(0)_{q}! =1$,
$${\binom{n}{i}}_{q} = \frac{(n)_{q}!}{(i)_{q}!(n - i)_{q}!}, \quad
\text{where }(n)_{q}! = \prod_{1 \le i \le n} (i)_{q}, \quad (n)_{q} =
\frac{q^{n} - 1}{q - 1}.$$

By the quantum binomial formula, if $ 1 \le n_j < N_j$, then $$\Delta
(x_{j}^{n_{j}}) = \sum_{0 \le i_{j} \le n_{j}}
{\binom{n_{j}}{i_{j}}}_{q_{j}}x_{j}^{i_{j}}\otimes  x_{j}^{n_{j}-i_{j}}.$$
We use the  following notation:  $$\bold n = (n_{1}, \dots, n_{j},
\dots , n_{\theta}), \quad  x^{\bold n}= x_{1}^{n_{1}} \dots
x_{j}^{n_{j}} \dots x_{\theta}^{n_{\theta}}, \quad \vert \bold n\vert
=  n_{1} + \dots  + n_{j} + \dots  + n_{\theta};$$  accordingly,
$\bold N= (N_{1}, \dots, N_{\theta})$, $\bold 1 = (1,\dots, 1)$. 
Also, we set $$\bold i \le \bold n \quad \text{if } i_{j} \le n_{j}, \, j
= 1, \dots, \theta.$$ Then, for  $\bold n \le\bold N - \bold 1$, we
deduce from the quantum binomial formula that  $$\Delta
(x^{\bold n}) = x^{\bold n}\otimes 1 + 1\otimes x^{\bold n} +
\sum\Sb 0 \le \bold i \le \bold n\\ 0 \ne \bold i \ne \bold n\endSb
c_{\bold n, \bold i}x^{\bold i}\otimes x^{\bold n-\bold i},$$ where
$c_{\bold n, \bold i} \ne 0$ for all $\bold i$. We shall need  
$$\aligned \Delta (x_{j}x_{i}) &= (x_{j} \otimes 1 + 1 \otimes
x_{j})(x_{i} \otimes 1 + 1 \otimes x_{i})\\ &=  x_{j}x_{i} \otimes 1
+ x_{j}\otimes x_{i} + \chi_{i}(g_{j}) x_{i}\otimes x_{j} + 1\otimes
x_{j}x_{i}, \endaligned \tag 3.1$$ for  $1 \le j, i \le \theta$.

\proclaim{Lemma 3.3} Let $R$ be a braided Hopf algebra in
${}_K^K\Cal{YD}$ satisfying hypothesis (A).  Then $\{x^{\bold n}:
\quad  \bold n \le \bold N - \bold 1\}$ is linearly independent.
Hence, $\dim R \ge N_{1}\dots N_{\theta}.$  In particular if any
element of $G(K)$ has order $p$, then $\dim R \ge p^{\theta}.$
\endproclaim

\demo{Proof} We shall prove by induction on $r$ that the set 
$$\{x^{\bold n}: \quad \vert \bold n\vert \le r, \quad \bold n \le
\bold N - \bold 1\}$$ is linearly independent.  

Let $r=1$ and let $a_{0} + \sum_{i=1}^{\theta} a_{i}x_{i} = 0$, wuth
$a_{j}\in k$, $0\le j \le \theta$. Applying $\varepsilon$, we see
that $a_{0} = 0$; by hypothesis we conclude that the other
$a_{j}$'s are also 0.

Now let $r> 1$ and  suppose that
$z = \sum_{\bold n:  \vert \bold n\vert \le r} a_{\bold n} x^{\bold
n} = 0$. Then $$0 = \Delta (z) =z\otimes 1 + 1\otimes z +
\sum_{ 1<\vert\bold n\vert \le r} a_{\bold n}\sum\Sb 0 \le \bold i
\le \bold n\\ 0 \ne \bold i \ne \bold n\endSb c_{\bold n, \bold
i}x^{\bold i}\otimes x^{\bold n-\bold i} = \sum_{1<\vert\bold n\vert
\le r} \,\sum\Sb 0 \le \bold i \le \bold n\\ 0 \ne \bold i \ne \bold
n\endSb a_{\bold n} c_{\bold n, \bold i}x^{\bold i}\otimes x^{\bold
n-\bold i} . $$ Now, if $\vert\bold n\vert \le r$, $0 \le \bold i \le
\bold n$, and $0 \ne \bold i \ne \bold n$, then $\vert\bold i\vert <
r$ and  $\vert\bold n - \bold i\vert < r$. By inductive hypothesis,
the elements $x^{\bold i}\otimes x^{\bold n-\bold i}$ are linearly
independent. Hence $a_{\bold n} c_{\bold n, \bold i} = 0$ and 
$a_{\bold n} = 0$ for all $\bold n$, $\vert \bold n \vert >1$. By the
step $r = 1$, $a_{\bold n} = 0$ for all $\bold n$.  \qed\enddemo

\bigpagebreak Let now $\theta \in \Bbb N$ and  $g_{1}, \dots,
g_{\theta} \in G(K)$, $\chi_{1}, \dots, \chi_{\theta} \in \Alg (A, k)$. 
We assume   that

\roster \item"(3.2)" the order $N_j$ of $q_{j}:= \chi_{j} (g_{j})$ is
finite. To avoid degenerate cases we also assume $N_{j} > 1$; {\it
cf.} Lemma 3.1.

\medpagebreak
 \item"(3.3)" $g_ig_j = g_jg_i$, \quad  $\chi_i\chi_j =
\chi_j\chi_i$, \quad for all $i,j$.

\medpagebreak
 \item"(3.4)" $\chi_i(h) g_i =  h_{(1)} \chi_i(h_{(2)}) g_i
\Ss(h_{(3)})$, for all $h\in K$, $1 \le i \le \theta$.

\medpagebreak
\item"(3.5)" $\chi_{j} (g_{i}) \chi_{i} (g_{j}) = 1$, \quad for all  $i
\ne j$. \endroster

\bigpagebreak

For $K = k(\Gamma)$ with $\Gamma$ finite abelian, the following
Lemma was  essentially proved in \cite{N, p. 1539}.

Let  $R$ be  the algebra   generated by  $x_{1},\dots, x_{\theta}$,
with relations  \roster \item"(3.6)" $x_{1}^{N_{1}} = 0, \dots,
x_{\theta}^{N_{\theta}} = 0,$

\medpagebreak
\item"(3.7)" $x_{i}x_{j} = \chi_{j}(g_{i}) x_{j}x_{i}$, if $i \ne j$.
\endroster

\proclaim{Lemma 3.4} $R$ has a unique braided Hopf algebra
structure in ${}_K^K\Cal{YD}$ such that the action and coaction are
determined by
 $$\delta (x_{j}) = g_{j}\otimes x_{j}, \qquad h. x_{j} =
\chi_{j}(h)x_{j} \quad \forall h\in \Gamma, \quad 1\le j \le \theta,$$ and
the $x_{i}$'s are primitive.  The dimension of $R$ is $N_{1}\dots
N_{\theta}.$ The coradical of $R$ is $k1$ and the space $P(R)$ of
primitive elements is the span of the $x_{i}'s$. $R$
is a coradically graded Hopf algebra, with respect
to the grading where the $x_{i}$'s are homogeneous
of degree 1.\endproclaim   
We denote this braided Hopf algebra by $\Cal R(g_1, \dots,
g_{\theta}; \chi_1, \dots, \chi_{\theta})$; it will be called  a {\it
quantum linear space\/} over $K$.

\demo{Proof} We first observe that $R$ is a $K$-module algebra
and  a $K$-comodule algebra because of conditions (3.3). Indeed,
we can extend the preceding action and coaction of $K$ to the
free algebra on generators $x_{1},\dots, x_{\theta}$; then we have
to see that the ideal generated by the   relations (3.6) and (3.7) is
stable by the action and coaction. This is clear  for (3.6); for (3.7),
it follows from (3.3). In addition, the Yetter-Drinfeld  condition on
$R$ holds because of, and indeed is equivalent to, (3.4).

We verify next that the elements $1\otimes x_{i} + x_{i}\otimes
1\in R\otimes R$ satisfy relations (3.6) and (3.7). The first follow
from the quantum binomial formula; the second, by direct
computation using (3.5). The counit is determined by 
$\varep(x_{i}) = 0$. The existence of the antipode follows from a
Lemma of Takeuchi, see \cite{M, 5.2.10}:  it is enough to check that
the the restriction of the identity to the coradical of $R$ is
invertible. But is not difficult to see that $R_0 = k1$. Indeed $R$ is
a graded  coalgebra whose homogeneous part of degree 0 is $k1$;
then use \cite{Sw, 11.1.1}. Thus $R$ is a braided Hopf algebra. By
Lemma 3.3, $\dim R \ge N_{1}\dots N_{\theta}$. But (3.5)
guarantees that the monomials  $\{x^{\bold n}: \quad  \bold n \le
\bold N - \bold 1\}$ generate $R$ as vector space; whence $\dim R
= N_{1}\dots N_{\theta}.$  

Finally, it is clear that $R = \oplus_{n \ge 0} R(n)$ is a graded Hopf
algebra, where $R(n)$ is generated by the monomials $x^{\bold n}$ such that
$\vert \bold n \vert = n$. Let $z\in P(R)$; we can assume that $z$ is
homogeneous. By the same argument as in the proof of Lemma 3.3, $n= 1$.
That is, $R_{1} = P(R)$. The last assertion follows from \cite{CM, 2.2}.
\qed\enddemo

\bigpagebreak  Quantum linear spaces are characterized by the
following Proposition.

\proclaim{Proposition 3.5} Let $R$ be a  braided Hopf algebra in
${}_K^K\Cal{YD}$ satisfying hypothesis (A).  Assume that $$\dim R
= N_{1}\dots N_{\theta}.$$  Then: \roster \item"(i)" $\chi_{j} (g_{i})
\chi_{i} (g_{j}) = 1$, for all $i\ne j$; and \item"(ii)" $R$ is a quantum
linear space. \endroster \endproclaim

\demo{Proof} Relations (3.6)  hold by Lemma 3.1. By Lemma 3.3,
$\{x^{\bold n}: \bold n \le \bold N - \bold 1\}$ is a basis of $R$,
which is then generated as an algebra by  $x_{1},\dots,
x_{\theta}$.  If $i >j$,  $x_{i}x_{j}$ can be expressed in the
following way: $$x_{i}x_{j} = \sum_{\bold n} c_{\bold n}x^{\bold
n},$$ for some $c_{\bold n} \in k$. Applying $\Delta$, we see that
$c_{\bold n} = 0$ unless $x^{\bold n} = x_{j}x_{i}$; so $x_{i}x_{j} = c
x_{j}x_{i}$, for some $c\in k$. By (3.1), we have $$  \multline
x_{i}x_{j} \otimes 1 + x_{i}\otimes x_{j} + \chi_{j}(g_{i})
x_{j}\otimes x_{i} + 1\otimes x_{i}x_{j} \\= cx_{j}x_{i} \otimes 1 +
cx_{j}\otimes x_{i} + c\chi_{i}(g_{j}) x_{i}\otimes x_{j} + 1\otimes
cx_{j}x_{i}.  \endmultline$$ By Lemma 3.3 again, $c =
\chi_{j}(g_{i})$ and $c\chi_{i}(g_{j}) = 1$. Hence (i) and relations
(3.7) hold. Applying the action and  coaction to both sides of the
equality (3.7), the conditions (3.3) follow.

We can define now a surjective  algebra homomorphism  $$ \Cal
R(g_1, \dots, g_{\theta}; \chi_1, \dots, \chi_{\theta}) @>>> R,$$
which is also a morphism of Yetter-Drinfeld modules. It is easy to
conclude that it is a homomorphism of braided Hopf algebras. By a
dimension argument, this map is an isomorphism. \qed \enddemo

\bigpagebreak
\subhead \S 4. Quantum linear spaces over abelian groups\endsubhead 

Let $\Gamma$ be a finite non-trivial  abelian group and let $H =
k(\Gamma)$.   We discuss in this Section the existence of quantum linear
spaces over $H$. 

Let  $\theta \in \Bbb N$.  A {\it datum for a quantum linear space} consists of
elements  $g_{1}, \dots, g_{\theta} \in \Gamma$, $\chi_{1}, \dots,
\chi_{\theta} \in \widehat{\Gamma}$ such that conditions (3.2), \dots, (3.5)
hold.  Explicitly, and because $\Gamma$ is abelian, we are then requiring the
following conditions:

\roster \item"(4.1)"  $q_{j}:= \chi_{j} (g_{j}) \neq 1$.

\medpagebreak
\item"(4.2)" $\chi_{j} (g_{i}) \chi_{i} (g_{j}) = 1$, \quad for all  $i
\ne j$. 
\endroster
We shall say that the datum, or its associated quantum linear space, has {\it
rank}  $\theta$. Given $\theta$, we are interested in describing all the data of rank $\theta$.
This description could be very cumbersome.
Let $\theta(\Gamma)$ be the greatest integer $\theta$ such that a datum of rank $\theta$ exists.

\proclaim{Lemma 4.1} Let $\Gamma = K\times H$, where $K$ and $H$ are
finite abelian groups. Then $\theta(\Gamma) \geq \theta(K) +
\theta(H)$. If the orders of $K$ and $H$ are coprime, then  $\theta(\Gamma) = \theta(K) + \theta(H)$.
\endproclaim

\demo{Proof} We identify $H$, $K$ with subgroups of $\Gamma$, and
$\widehat H$, $\widehat K$ with subgroups of $\widehat{\Gamma}$. Let
$h_{1}, \dots, h_{\mu}$, $\eta_{1}, \dots, \eta_{\mu}$ be a datum for $H$ and
let $k_{1}, \dots, k_{\nu}$, $\zeta_{1}, \dots, \zeta_{\nu}$ be a datum for
$H$. Then  $$h_{1}, \dots, h_{\mu}, k_{1}, \dots, k_{\nu} \text{ in }
\Gamma, \qquad \eta_{1}, \dots, \eta_{\mu}, \zeta_{1}, \dots, \zeta_{\nu}
 \text{ in } \widehat{\Gamma},$$
is clearly a datum for $\Gamma$. Hence $\theta(\Gamma) \geq \theta(K) +
\theta(H)$. 

Conversely, asume that the orders of $H$ and $K$ are coprime and let $g_{1},
\dots, g_{\theta} \in \Gamma$, $\chi_{1}, \dots, \chi_{\theta} \in
\widehat{\Gamma}$ be a datum for $\Gamma$. Let us decompose 
$$
g_{i} = h_{i}k_{i}, \text{ where } h_{i} \in H, k_{i} \in K, \text{ and }
\chi_{i} = \eta_{i}\zeta_{i}, \text{ where } \eta_{i} \in \widehat H, 
\zeta_{i} \in \widehat K, \quad 1 \le i \le \theta. $$
We claim that $h_{i}$, $\eta_{i}$, $i\in I$, where   $I :=\{i: \,\eta_{i}(h_{i})
\neq 1\}$ is a datum for $H$, and,  similarly, that $k_{i}$, $\zeta_{i}$, $i\in
J$, where   $J :=\{i: /, \zeta_{i}(k_{i}) \neq 1\}$ is a datum for $K$. Clearly,
$$\chi_{i} (g_{i}) \neq 1 \text{  implies  } \eta_{i}(h_{i})
\neq 1\text{ or } \zeta_{i}(k_{i}) \neq 1;$$
that is, the claim implies $\theta(\Gamma) \leq \theta(K) + \theta(H)$.
We check then the claim. Condition (4.1) is forced by the choice of the index
sets. For (4.2), observe that
$$
1 = \chi_{j} (g_{i}) \chi_{i} (g_{j}) = \eta_{j} (h_{i}) \eta_{i}(h_{j})
\zeta_{j} (k_{i}) \zeta_{i} (k_{j}) $$
implies $1 = \eta_{j} (h_{i}) \eta_{i}(h_{j}) = \zeta_{j} (k_{i}) \zeta_{i}
(k_{j})$, because the orders of $\eta_{j} (h_{i}) \eta_{i}(h_{j})$ and
$\zeta_{j} (k_{i}) \zeta_{i} (k_{j})$ are coprime. \qed\enddemo

By the preceding Lemma, we are reduced to investigate the behaviour of
$\theta(\Gamma)$ when $\Gamma$ is an abelian $p$-group, $p$ a prime.

\proclaim{Lemma 4.2}Let $\Gamma$ be a cyclic $p$-group, where $p$ is an
odd prime. Then $\theta(\Gamma) = 2$. \endproclaim
\demo{Proof} We first prove that $\theta(\Gamma) \leq 2$. It is enough to
show that no datum of rank 3 exists.  Let us assume, on the contrary, that
$g_{1}, g_{2}, g_{3}\in \Gamma$, $\chi_{1}, \chi_{2}, \chi_{3}\in
\widehat{\Gamma}$, satisfy (4.1), (4.2). Let $g$ be a generator of the
subgroup $\langle g_{1}, g_{2}, g_{3}\rangle$ and let $p^{s}$ be the
order of $g$, where $s$ is a positive integer. Let $\zeta$ be a primitive
$s$-th root of 1. Let $a_{1}, a_{2}, a_{3}$, $b_{1}, b_{2}, b_{3}$ be integers
such that
$$g_{i} = g^{b_{i}}, \quad \chi_{i}(g) = \zeta^{a_{i}}, \quad 1\le i \le 3.$$
Then condition (4.1) means that $a_{i}b_{i} \not\equiv 0 \mod p^{s}$ and (4.2)
that
$$
\align
a_{1}b_{2} + a_{2}b_{1} &\equiv 0 \mod p^{s}\tag 4.3\\
a_{1}b_{3} + a_{3}b_{1} &\equiv 0 \mod p^{s}\tag 4.4\\
a_{2}b_{3} + a_{3}b_{2} &\equiv 0 \mod p^{s}.\tag 4.5
\endalign
$$
On the other hand, there exist integers $r_{1}, r_{2}, r_{3}$ such that 
$$b_{1}r_{1} + b_{2}r_{2} + b_{3}r_{3} \equiv 1 \mod p^{s}. \tag 4.6$$
Now, we multiply (4.3) by $b_{3}$, (4.4) by $b_{2}$, (4.5) by $b_{1}$ and
conclude (since $p$ is odd) that 
$$a_{1}b_{2}b_{3}\equiv 0 \mod p^{s}, \quad a_{2}b_{1}b_{3}\equiv 0 \mod
p^{s}, \quad a_{3}b_{1}b_{2}\equiv 0 \mod p^{s}.$$ 
Let us write
$$a_{1} = p^{t}\widetilde{a_{1}}, \quad \text{ where } t\geq 0, p\not \vert
\widetilde{a_{1}}.$$
Then $p^{s-t} \vert b_{2}b_{3}$ and hence there exist positive integers $h$,
$j$ such that $p^{h} \vert b_{2}$, $p^{j} \vert b_{3}$ and $h + j = s-t$. From
(4.3), (4.4) we deduce that $p^{t+h} \vert a_{2}b_{1}$ and $p^{t+j} \vert
a_{3}b_{1}$. Now assume that $p\not \vert b_{1}$. Then
$$
p^{t+h} \vert a_{2} \implies p^{t+ 2h} \vert a_{2}b_{2} \implies t+ 2h < s,$$
and similarly, $t+ 2j < s$. But then $h < \dfrac{ s - t}{2}$, $j < \dfrac{ s -
t}{2}$ and therefore $h + j < s$, which is not possible. Hence $p \vert b_{1}$.
By symmetry, $p \vert b_{2}$, $p \vert b_{3}$. This contradicts (4.6) and
finishes the proof of $\theta(\Gamma) \leq 2$.

Let $g$ denote now a generator of $\Gamma$  and let again $p^{s}$ be the
order of $g$ and $\zeta$  a primitive $s$-th root of 1. Then $g_{1} = g$,
$g_{2} = g^{a}$, $\chi_{1}$ given by $\chi_{1}(g) = \zeta$ and $\chi_{2}$
given by $\chi_{2}(g) = \zeta^{-a}$ is a datum of rank 2 whenever $a^{2}
\not\equiv 0 \mod p^{s}$. \qed\enddemo

\bigpagebreak
\subhead \S 5. Pointed Hopf algebras whose diagrams are quantum linear
spaces \endsubhead

Let $\Gamma$ be a finite abelian group. We fix a decomposition $\Gamma =
\langle  y_{1}\rangle\oplus \dots\oplus  \langle  y_{\sigma}\rangle$ and
we denote by $M_{\ell}$ the order of $y_{\ell}$, $1\leq \ell \leq \sigma$.

\medpagebreak
Let $g_{1}, \dots, g_{\theta} \in \Gamma$, $\chi_{1}, \dots, \chi_{\theta}
\in \widehat{\Gamma}$ be a datum for quantum linear space; {\it i. e.},  (4.1),
(4.2) hold. We set $q_{i} = \chi_{i}(g_{i})$, $N_{i}$  the order of $q_{i}$. We
abbreviate $\Cal R :=  \Cal R(g_1, \dots, g_{\theta}; \chi_1, \dots,
\chi_{\theta})$ for the quantum linear space defined in Section 3.

\medpagebreak
A compatible datum $\Cal D$ for $\Gamma$ and $\Cal R$ consists of
\roster   
\item"(5.1)" a scalar $\mu_{i} \in \{0, 1\}$ for each $i$, $1\leq i \leq
\theta$;  it is arbitrary if $g_{i}^{N_{i}} \neq 1$ and $\chi_{i}^{N_{i}} = 1$,
but 0 otherwise;
\item"(5.2)" a scalar $\lambda_{ij} \in k$ for each $i, j$, $1\leq i < j \leq
\theta$;  it is is arbitrary if $g_{i}g_{j} \neq 1$ and $\chi_{i}\chi_{j} = 1$,
but 0 otherwise. \endroster

\remark{Remark}   If $\lambda_{ij}\neq 0$ and  $\lambda_{ih}\neq 0$,
then $\chi_{i}\chi_{j} = 1$ and $\chi_{i}\chi_{h} = 1$; hence $\chi_{j} =
\chi_{h}.$  If in addition the order $N_i$ of $q_{i}:=\chi_{i} (g_{i})$ is odd,
then  $j = h$. Indeed, suppose $j \neq h$. Then
$$1 = \chi_{j}(g_{h})\chi_{h}(g_{j}) = \chi_{i}(g_{h})^{-1}\chi_{i}(g_{j})^{-1}
=  \chi_{i}(g_{i})^{-2}.$$
Here, the first equality is by (4.2); the second, because $\chi_{j} = \chi_{h}
= \chi_{i}^{-1}$; the third, because  $\chi_{i}(g_{i})^{-1} = \chi_{j}(g_{i}) =
\chi_{i}(g_{j})^{-1}$ and similar with $h$ instead of $j$. Now $N_{i}$ odd
forces $1 =  \chi_{i}(g_{i})$, which is excluded by  (3.2). 
\endremark

\bigpagebreak
Let $\eta$ be the  injective map  from $\Gamma$ to
the free algebra $k\langle h_{1},\dots, h_{\sigma}, a_{1}, \dots, a_{\theta} 
\rangle$ given by 
$$\eta(y_{1}^{n_{1}}, \dots, y_{\sigma}^{n_{\sigma}}) = h_{1}^{n_{1}} \dots
h_{\sigma}^{n_{\sigma}}, \quad 0\leq n_{\ell}\leq M_{\ell} - 1, \quad 1\leq
\ell \leq \sigma.$$
We shall identify elements of  $\Gamma$ with elements of the free algebra
$k\langle h_{1},\dots, h_{\sigma}, a_{1}, \dots, a_{\theta}  \rangle$ via
$\eta$ without further notice.

\definition{Definition} Let $\Gamma$ be a finite abelian group, $\Cal R$ a
quantum linear space and $\Cal D$ a compatible datum. Keep the notation
above. Let $\Aq$ be the algebra presented by generators $h_{\ell}$, $1\leq
\ell \leq \sigma$, and $a_{i}$, $1\leq i \leq \theta$ with defining relations

\medpagebreak
\roster \item"(5.3)" $h_{\ell}^{M_{\ell}} = 1$, $1\leq \ell \leq \sigma$; 

\medpagebreak
\item"(5.4)" $h_{\ell}h_{t} = h_{t}h_{\ell}$, $1\leq t <\ell \leq \sigma$;

\medpagebreak
\item"(5.5)" $a_{i}h_{\ell} = \chi_{i}^{-1}(y_{\ell})h_{\ell} a_{i}$, $1\leq
\ell \leq \sigma$, $1\leq i \leq \theta$;

\medpagebreak
\item"(5.6)" $a_{i}^{N_{i}} = \mu_{i}\left(1 -
g_{i}^{N_{i}}\right)$, $1\leq i \leq \theta$;

\medpagebreak
\item"(5.7)" $a_{j}a_{i} = \chi_{i}(g_{j}) a_{i}a_{j} + \lambda_{ij} \left(1 -
g_{i} g_{j} \right)$, $1\leq i < j \leq \theta$. \endroster
\enddefinition

We shall denote in  what follows by the same letters the generators of the
free algebra $k\langle h_{1},\dots, h_{\sigma}, a_{1}, \dots, a_{\theta} 
\rangle$ and their classes in $\Aq$; no trouble should arise.

\proclaim{Lemma 5.1} There exists a unique Hopf algebra structure on 
$\Aq$ such that  $$\Delta(h_{\ell}) = h_{\ell}\otimes h_{\ell},
\qquad \Delta(a_{i}) = a_{i} \otimes 1 + g_{i}\otimes a_{i}, \quad 1\leq \ell
\leq \sigma,\quad 1\leq i \leq \theta. \tag 5.8$$ \endproclaim 
\demo{Proof} Let also $\Delta: k\langle h_{1},\dots, h_{\sigma}, a_{1}, \dots, a_{\theta} 
\rangle \to \Aq\otimes \Aq$ denote the algebra map defined by (5.8). 
Clearly, $\Delta(h) = h\otimes h$ whenever $h$ is a monomial in the
$h_{\ell}$'s.  We have to verify that the
elements $H_{\ell} = h_{\ell}\otimes h_{\ell}$, $A_{i} = a_{i} \otimes 1 +
g_{i}\otimes a_{i}$ satisfy the defining relations.  This is not difficult for
(5.3), (5.4), (5.5).  For relations (5.6), (5.7), the reason is the same: both
sides of each equality are skew-primitive elements related to the same
group-likes. For instance, we have  $$\Delta(a_{i})^{N_{i}} = a_{i}^{N_{i}}
\otimes 1 + g_{i}^{N_{i}}\otimes a_{i}^{N_{i}} = \mu_{i}(1 - g_{i}^{N_{i}}) \otimes 1 +
g_{i}^{N_{i}} \otimes \mu_{i}(1 - g_{i}^{N_{i}}) = \Delta\mu_{i}(1 -
g_{i}^{N_{i}}).$$ Here the first equality follows from  (5.4) and the
definition of $\Delta$ via the quantum binomial formula, since the order of
$q_{i}$ is $N_{i}$; the second, from (5.6); the third is clear. This proves that
$H_{\ell}$, $A_{i}$ satisfy (5.6). For (5.7), the computation is also direct. It is
clear that $\Delta$ is coassociative.

\medpagebreak
The algebra map $\varep: \Aq \to k$ uniquely determined
by $\varep(h_{\ell}) = 1$, $\varep(a_{i}) = 0$, for all $\ell$ and $i$, is the
counit of $\goth \Aq$.   We claim that there is a unique algebra map $\Ss:
\goth \Aq \to \Aq^{op}$ such that  for all $\ell$ and $i$,
$$\Ss(h_{\ell}) = h_{\ell}^{-1}, \qquad \Ss(a_{i}) = -g_{i}^{-1}a_{i}.$$
The verification of relations (5.3), (5.4), (5.5), (5.7) is straightforward.
For (5.6), we first check by induction that
$$\Ss (a_{i})^{n} = (-1)^{n}q_{i}^{n(n-1)/2} g_{i}^{-n}a_{i}^{n}.$$
As $q_{i}$ is a primitive $N_{i}$-th root of 1, $(-1)^{N_{i}}q_{i}^{N_{i}(N_{i}-1)/2}  =-1$. Hence
$$\Ss (a_{i})^{N_{i}} = - g_{i}^{-N_{i}}a_{i}^{N_{i}} = \mu_{i}(1 - g_{i}^{-N_{i}}) = \mu_{i}(1 -
\Ss(g_{i})^{N_{i}}). $$ 
The map $\Ss$ is clearly an antipode and the Lemma follows. \qed\enddemo

\proclaim{Proposition 5.2}  Let $\Gamma$ be a finite abelian group, $\Cal R$
a quantum linear space and $\Cal D$ a compatible datum. Keep the notation
above. The set of monomials
$$h_{1}^{r_{1}}\dots h_{\sigma}^{r_{\sigma}} a_{1}^{s_{1}}\dots
a_{\theta}^{r_{\theta}}, \qquad 0\le r_{\ell} < M_{\ell}, \quad  0\le s_{i}
< N_{i}, \quad 1\leq \ell \leq \sigma, \quad 1\leq i \leq \theta $$ is a basis
of $\Aq$. In particular,  $$\dim \Aq = \prod_{1\leq \ell \leq \sigma}
M_{\ell} \prod_{1\leq i \leq \theta} N_{i} = \vert \Gamma\vert \dim \Cal R.
\tag 5.9$$ \endproclaim \demo{Proof} Let us assume that the scalars
$\mu_{i}$, $\lambda_{ij}$ are arbitrary. It is not difficult to conclude from
relations (5.4), (5.5) and (5.7) that these monomials generate the vector
space $\Aq$. By the Diamond Lemma \cite{Be}, it is then enough to verify
that the following overlaps can be reduced to the same normal form: 

\roster \item"(5.10)"  \quad $\left(a_{i}h_{\ell}\right)h_{\ell}^{M_{\ell} -
1} = a_{i} \left(h_{\ell}h_{\ell}^{M_{\ell} - 1}\right)$;
\item"(5.11)" \quad $\left(a_{i}h_{\ell}\right)h_{t} = a_{i} \left(h_{\ell}
h_{t} \right)$, $t < \ell$;
\item"(5.12)" \quad $\left(a_{i}^{N_{i} - 1}a_{i}\right)h_{\ell} = a_{i}^{N_{i}
- 1} \left(a_{i} h_{\ell} \right)$;
\item"(5.13)"  \quad $\left(a_{j}^{N_{j} - 1}a_{j}\right)a_{i} = a_{j}^{N_{j}
- 1} \left(a_{j} a_{i} \right)$, $i<j$;
\item"(5.14)" \quad $\left(a_{j}a_{i}\right)a_{i}^{N_{i} - 1} = a_{j}
\left(a_{i} a_{i}^{N_{i} - 1} \right)$, $i<j$; 
\item"(5.15)" \quad $\left(a_{j}a_{i}\right)h_{\ell} = a_{j}
\left(a_{i} h_{\ell} \right)$, $i<j$.
\endroster
Here we order the monomials in the following way. If $z_{1}<z_{2}< \dots <
z_{m}$ are indeterminates, we define the standard ordering on monomials
$A, B$ in $z_{1},\dots,z_{m}$ in the usual way: $A<B$ if
$\length(A)<\length(B)$, or $A$ and $B$ have the same length and $A$ is
lexicographically smaller than $B$.  If $A$ is a monomial in $h_{1}, \dots ,
a_{\theta}$, let $\phi(A)$ be its $a$-part, that is the image under the
monoid homomorphism $\phi$ with $\phi(h_{\ell}) = 1$, $\phi(a_{i}) = a_{i}$
for all $\ell$, $i$.  We order the monomials in $h_{1},\dots,a_{\theta}$ as
follows: $h_{1} < \dots < h_{\sigma} < a_{1} < \dots <a_{\theta}$; $A<B$ if
$\phi(A)<\phi(B)$ in the standard ordering of the monomials in $a_{1}, \dots,
a_{\theta}$, or $\phi(A) = \phi(B)$ and $A$ is smaller than $B$ in the
standard ordering of $h_{1}, \dots, a_{\theta}$.

The verification of (5.10), (5.11) is easy and gives no condition.  The
verification of (5.12) amounts to 
$$ \mu_{i}\left(1 - g_{i}^{N_{i}}\right)h_{\ell}  = \mu_{i}
\chi_{i}^{-1}(y_{\ell})^{N_{i}} h_{\ell} \left(1 - g_{i}^{N_{i}}\right). $$
This imposes the condition
$$\text{ If } g_{i}^{N_{i}} \neq 1 \text{ and } \chi_{i}^{N_{i}} \neq 1
\text{ then } \mu_{i} = 0.\tag 5.16$$ 
The verification of (5.13) turns to 
$$\multline \mu_{j}\left(1 - g_{j}^{N_{j}}\right) a_{i} \\ = 
\mu_{j}\chi_{i}( g_{j})^{N_{j}} a_{i} - \mu_{j} g_{j}^{N_{j}} a_{i} +
\lambda_{ij}\left(1 + \chi_{i}( g_{j}) + \chi_{i}( g_{j})^{2} + \dots
+ \chi_{i}( g_{j})^{N_{j} - 1}\right) a_{j}^{N_{j} - 1}; \endmultline$$
and so we need the conditions
$$\align
\lambda_{ij}\left(1 + \chi_{i}( g_{j}) + \chi_{i}( g_{j})^{2} + \dots
+ \chi_{i}( g_{j})^{N_{j} - 1}\right) &= 0. \tag 5.17 \\
\text{ If } g_{j}^{N_{j}} \neq 1 \text{ and } \chi_{i}(g_{j})^{N_{j}} \neq 1
\text{ then } \mu_{i} &= 0.\tag 5.18\endalign$$
In the same vein, for (5.14) and (5.15) it is necessary that
$$\align
\lambda_{ij}\left(1 + \chi_{i}( g_{j}) + \chi_{i}( g_{j})^{2} + \dots
+ \chi_{i}( g_{j})^{N_{i} - 1}\right) &= 0. \tag 5.19 \\
\text{ If } g_{i}^{N_{i}} \neq 1 \text{ and } \chi_{i}(g_{j})^{N_{i}} \neq 1
\text{ then } \mu_{i} &= 0.\tag 5.20\\
\text{ If } g_{i}g_{j} \neq 1 \text{ and } \chi_{i}\chi_{j} \neq 1
\text{ then } \lambda_{ij} &= 0.\tag 5.21\endalign$$
Now it is harmless to assume that 
$$\mu_{i} = 0 \text{ if }g_{i}^{N_{i}} = 1, \lambda_{ij}  = 0 \text{ if } g_{i}g_{j}
=1.$$
 The combination of  this last assumption and (5.16)  is exactly the
constraint in (5.1); in turn, the constraint in (5.2) is equivalent to the
assumption together with (5.21).  Also, condition (5.16) implies (5.18) and
(5.20). It remains to show that (5.17) and (5.19) are consequences of (5.21). 

Indeed, assume that $\lambda_{ij} \neq 0$; by (5.21), this is only possible if
$g_{i}g_{j} \neq 1$ and $\chi_{i}\chi_{j} = 1$. But then $\chi_{i}(g_{i}) =
\chi_{j}(g_{j})^{-1}$, thanks to (4.2). Thus $N_{i} = N_{j}$. Moreover, 
$\chi_{i}(g_{j})^{N_{j}} \chi_{j}(g_{j})^{N_{j}} = 1$ and hence
$\chi_{i}(g_{j})^{N_{j}}  = 1$. Therefore (5.17) and (5.19) hold. \qed\enddemo

\bigpagebreak
\proclaim{Corollary 5.3} The Hopf algebra $\Aq$ is pointed and its
coradical filtration is given by
$$\Aq_{n} = \left\langle h_{1}^{r_{1}}\dots h_{\sigma}^{r_{\sigma}}
a_{1}^{s_{1}}\dots a_{\theta}^{s_{\theta}}, \quad 0\le r_{\ell} < M_{\ell},
\quad  0\le s_{i} < N_{i}, \quad \forall \ell, i,  \;  \sum_{i}
s_{i} \leq n\right\rangle. \tag 5.22$$ In particular,  $$\aligned 
P_{g_{i}, 1} (\Aq) &= \left( \oplus_{j: g_{j} = g_{i}}k(1-g_{j})\right) \oplus \left( \oplus_{j: g_{j} = g_{i}} ka_{j}\right), \quad 1\leq i \leq \theta, \\
P_{g, 1} (\Aq) &= k(1-g) \quad \text{ if } g\neq g_{i}.
\endaligned \tag 5.23$$

\endproclaim

\demo{Proof}The subalgebra $k[h_{1},\dots, h_{\sigma}]$ of $\Aq$ coincides
with its coradical. Indeed, $k[h_{1},\dots, h_{\sigma}] \supset
\Aq_{0}$ by \cite{M, 5.5.1} and the other inclusion is evident. Hence, $\Aq$
is pointed and  $\Aq_{0}$ is isomorphic to the group algebra of $\Gamma$.

Now we consider the graded Hopf algebra $\gr \Aq$ associated to the
coradical filtration, and the diagram  of $\Aq$. It follows from Proposition
5.2 that  the diagram is isomorphic to $\Cal R$. By Lemma 3.4, we know the
coradical filtration of $\Cal R$. By Lemmas 2.3 and 2.4, we  know  then the
coradical filtration of $\gr \Aq$. We conclude, by a recursive argument that
the coradical filtration of $\Aq$ is given by (5.22).  In particular, 
$$\Aq_{1} = \Aq_{0} \oplus \Aq_{0}a_{1} \oplus \Aq_{0}a_{2} \dots \oplus
\Aq_{0}a_{\theta}.$$ The claim (5.23) follows by a direct computation.
\qed\enddemo

\bigpagebreak
Let now $A$ be a finite dimensional  pointed Hopf algebra such that the
group $G(A)$ of its group-like elements is isomorphic to $\Gamma$.  We
denote $H = k(\Gamma)$. By the Theorem of Taft and Wilson \cite{M, Thm.
5.4.1}, $A_1 = k(\Gamma) + (\oplus_{g,h \in \Gamma} P_{g,h})$.  

If $M$ is an $H$-module (respectively, comodule) then $M^{\chi}$ (resp.,
$M_{g}$) denotes the isotypic component of type $\chi\in
\widehat{\Gamma}$ (resp., of type $g\in \Gamma$). If $M$ is an object in  
${}_H^H\Cal{YD}$ then $M_{g}^{\chi} := M_{g} \cap M^{\chi}$. Any finite
dimensional $M \in {}_H^H\Cal{YD}$ decomposes as $$M = \oplus_{g\in
\Gamma, \chi\in \widehat{\Gamma}}M_{g}^{\chi}.$$

The adjoint action of $\Gamma$ on $A$ leaves stable
each space $P_{g,h}$; hence, we can further decompose  $P_{g,h} =
\oplus_{\chi \in \widehat{\Gamma}}  P_{g,h}^{\chi}$.

\proclaim{Lemma 5.4} Let $\gr A$ be the graded Hopf algebra associated to
the coradical filtration and let $R$ be the diagram of $A$.

(i). The first term of the coradical filtration of
$A$ is given by $$A_1 = k(\Gamma) \oplus (\oplus\Sb g,h\in
\Gamma \\ \chi \in \widehat{\Gamma} , \chi\neq \epsilon \endSb
P_{g,h}^{\chi}).$$Thus the second summand is isomorphic to $\gr
A(1)$. 

\medpagebreak
(ii). If $P(R) = \oplus_{1\le j \le M} P(R)^{g_{i}}$ with $P(R)^{g_{i}}
\neq 0$, then $P_{g,h}(A)$ contains properly $P_{g,h}(A)\cap
k(\Gamma) = k(g-h)$ if and only if $(g, h) = (g_{i}s, s)$, for some
$s\in \Gamma$.\endproclaim    \demo{Proof} If $\epsilon$ is the
trivial character of $\Gamma$, then $P_{g,h}^{\epsilon}\subset
k(\Gamma)$ by Lemma 3.1. Since $\Gamma$ is abelian,
$P_{g,h}^{\epsilon} = P_{g,h} \cap k(\Gamma)$. This shows part (i). 
Part (ii) follows at once from part (i) and formulas (2.1).
\qed\enddemo

\proclaim{Lifting Theorem 5.5} Let $A$ be a pointed finite dimensional
Hopf algebra with coradical $H = k(\Gamma)$, where $\Gamma$ is an abelian
group as above. Let $\gr A$ be the graded Hopf algebra associated to the
coradical filtration. Let $\Cal R$ be the diagram of $A$. We
assume that $\Cal R$ is a quantum linear space.  

Then there exists a compatible datum $\Cal D$ such that $A$ is isomorphic to
$\Aq$ as Hopf algebras. \endproclaim

\demo{Proof}Let $x_{1},\dots, x_{\theta}$ be the generators of $\Cal R$
satisfying the relations (3.6), (3.7). We identify $x_{j}$, resp. $h\in
\Gamma$, with $x_{j} \# 1$, resp. $1\# h$, in $\Cal R\# k(\Gamma)\simeq
\gr A$. By (2.1), we see that $\gr A$ can be presented by
generators  $h_{\ell}$, $1\leq \ell \leq \sigma$, $x_{i}$,   $1\le i \le
\theta$ and relations (5.3),
(5.4), \roster \item"(5.24)" $x_{i}^{N_{i}} = 0$,
\item"(5.25)" $ h_{\ell} x_{i} =\chi_{i}(h_{\ell})x_{i}h_{\ell}$,
\item"(5.26)" $x_{i}x_{j} - \chi_{j}(g_{i})x_{j}x_{i} = 0$,
\endroster  
for all $1\leq \ell \leq \sigma$, $1\le i \neq j \le \theta$. The Hopf algebra
structure of $\gr A$ is determined by  $$\Delta(h) = h\otimes h, \qquad
\Delta(x_{i}) = x_{i}\otimes 1 + g_{i}\otimes x_{i},$$ $1\le i \le
\theta$, $h\in \Gamma$. Hence $x_{i} \in P_{g_{i},1}(\gr
A)^{\chi_{i}}$. According to Lemma 5.4, we can choose $a_{i} \in
P_{g_{i},1}(A)^{\chi_{i}}$ such that  $\overline{a_{i}} = x_{i}$ in $\gr A (1) =
A_{1}/A_{0}$. By Lemma 2.2, $A$ is generated by $h_{\ell}$, $1\leq \ell \leq
\sigma$ $a_{i}$,  $1\le i \le \theta$.  It is clear that relations (5.3) and
(5.4) also hold in $A$. We verify now that  relations (5.5), (5.6), (5.7) hold
for some collection of scalars $\mu_{i}$, $\lambda_{ij}$, and at the same
time, that this choice must fulfill the constraints in (5.1) and (5.2). For
(5.5), this follows from  the choice of the $a_{i}$'s.  We check (5.6). By the
quantum binomial formula,  $$a_{i}^{N_{i}} \in
P_{g_{i}^{N_{i}},1}(A)^{\chi_{i}^{N_{i}}}.$$ We know that
$$g_{i}a_{i}^{N_{i}}g_{i}^{-1} = \chi_{i}^{N_{i}}(g_{i}) a_{i}^{N_{i}} =
q_{i}^{N_{i}}a_{i}^{N_{i}} = a_{i}^{N_{i}};$$ by Lemma 3.1,  $a_{i}^{N_{i}} \in
k(g_{i}^{N_{i}} - 1)$. Dividing out $a_{i}$ by an appropiate scalar, we see
that relations (5.6) hold, for  $\mu_{i}$ either 0 or 1. If $g_{i}^{N_{i}} \neq
1$ we can assume without trouble that $\mu_{i} = 0$. So let us suppose
that $g_{i}^{N_{i}} = 1$. If $\mu_{i} = 1$ then 
$$h_{\ell}a_{i}^{N_{i}}h_{\ell}^{-1} = \chi_{i}^{N_{i}}(h_{\ell}) a_{i}^{N_{i}} =
 a_{i}^{N_{i}};$$
hence $\chi_{i}^{N_{i}}$ is forced to be 1.

\medpagebreak
We prove now (5.7). By (4.2) and the choice of the $a_{i}$'s, it follows that
$$a_{i}a_{j} - \chi_{j}(g_{i})a_{j}a_{i}  \in P_{1,
g_{i}g_{j}}(A)^{\chi_{i}\chi_{j}}.$$ 
By Lemma  5.4, if $a_{i}a_{j} - \chi_{j}(g_{i})a_{j}a_{i} \notin
k(\Gamma)$, then for some $h\neq i, j$, $\chi_{i}\chi_{j} = \chi_{h}$. But
then  $g_{i}g_{j} = g_{h}$. By (4.2) again, 
$$1 = \chi_{h}(g_{i}) \chi_{i}(g_{h}) = \chi_{i}(g_{i}) \chi_{j}(g_{i})
\chi_{i}(g_{i}) \chi_{i}(g_{j}) = \chi_{i}(g_{i})^{2}$$ 
and hence $\chi_{i}(g_{i}) = -1$. Similarly, $\chi_{j}(g_{j}) = -1$. So
$$\chi_{h}(g_{h}) = \chi_{i}(g_{i}) \chi_{j}(g_{i})
\chi_{i}(g_{j}) \chi_{j}(g_{j}) = 1,$$
a contradiction. Therefore
$a_{i}a_{j} - \chi_{j}(g_{i})a_{j}a_{i} \in k(\Gamma)$ and by Lemma 3.1,
there exist scalars $\lambda_{ij}$ such that $a_{i}a_{j} - \chi_{j}(g_{i})a_{j}a_{i}
= \lambda_{ij}\left(1 - g_{i}g_{j}\right)$; {\it i.e.} (5.7) holds. If $g_{i}g_{j}
= 1$ we assume without harm that $\lambda_{ij} = 0$. If $g_{i}g_{j} \neq 1$
and $\lambda_{ij} \neq 0$ then, arguing as for the
$\mu_{i}$'s, we see that $\chi_{i}\chi_{j} = 1$. Hence the collection
$\lambda_{ij}$ satisfies the constraints of (5.2). 

Then the datum $\Cal D = \left(\mu_{i}, \lambda_{ij}\right)$ is compatible
and  we have  a Hopf algebra surjection $\Aq \to A$. As $\Aq$ and $A$
have the same dimension, they are isomorphic. \qed\enddemo

We deduce now Theorem 0.2 from Theorem 5.5. We need the following
Lemma.

\proclaim{Lemma 5.6}Let $\Gamma$ be a finite non-trivial  abelian group
and let $H = k(\Gamma)$.  Let $R$ be a  braided Hopf  algebra in
${}_H^H\Cal{YD}$, with trivial coradical: $R_{0} = R(0) = k1$. 

\roster \item"(a)" If $\dim R = p$ then $\dim P(R) = 1$ and $R$ is a
quantum line. \item"(b)" If $\dim R = p^{2}$ then $\dim P(R) = 1$ or
2, and $R$ is respectively a quantum line or a quantum plane.
\endroster \endproclaim \demo{Proof} Let $R$ be a finite
dimensional braided Hopf  algebra in ${}_H^H\Cal{YD}$, with trivial
coradical.  Since $R_{0} = k1$ and $R\supsetneq R_{0}$, $P(R) \neq
0$. On the other hand, $P(R)$ is a Yetter-Drinfeld submodule of $R$,
hence $P(R)  = \oplus_{g\in \Gamma, \chi\in
\widehat{\Gamma}}P(R)_{g}^{\chi}.$  

Let   $x \in P(R)_{g}^{\chi}$, $x\neq 0$, for some $g\in \Gamma$,
$\chi\in \widehat{\Gamma}$. Let $q = \chi(g)$ and let $N$ be the
order of $q$;  $q \ne 1$ by Lemma 3.1; that is, $N >1$. It is not
difficult to see that the subalgebra $k[x]$ of $R$ is a braided Hopf
subalgebra of dimension $N$.  It follows from the Nichols-Zoeller
Theorem that $N$ divides the dimension of $R$, see \cite{AS2,
Proposition 4.9}. 

Let $x_{1}, \dots, x_{\theta}$ be a basis of $P(R)$ such that $x_{j}
\in P(R)_{g_{j}}^{\chi_{j}}$, for some $g_{j}\in \Gamma$, $\chi_{j}\in
\widehat{\Gamma}$, for all $j$. Let $N_{j}$ be the order of
$\chi_{j}(g_{j})$.

If the dimension of $R$ is $p$, the considerations above show that
$R = k[x_{1}]$. This proves part (a). 

We now assume that the dimension of $R$ is $p^{2}$. If $N_{1} =
p^{2}$, then $\theta=1$ and $R$ is a quantum line. So we can
further suppose that $N_{j} = p$ for all $j$. By Lemma 3.3,
$\theta\leq 2$. If $\theta = 1$, then Theorem 3.2 forces $\dim R =
p$. This is a contradiction and therefore $\theta = 2$. We conclude
then, by Proposition 3.5, that $R$ is a quantum plane.
\qed\enddemo

\demo{Proof of Theorem 0.2} Let $\gr A$ be the graded Hopf algebra
associated to the coradical filtration and let $R$ be the braided Hopf
algebra in ${}_H^H\Cal{YD}$ such that $\gr A \simeq R\# H$ as in {\it 2.2.} If
the index of $H$ in $A$ is $p$ or $p^{2}$, then $R$ is a quantum line or plane,
according to Lemma 5.6. The description follows now from Theorem 5.5. 
\qed\enddemo

\bigpagebreak
\subhead \S 6. Families of Hopf algebras of the same dimension\endsubhead 
We shall specialize Proposition 5.2 to the simplest possible $\Gamma$ and
$\Cal R$ and suitable $\Cal D$.

\bigpagebreak
Let us  assume that $\Gamma$ is a cyclic group of order $MN$, where
$M > 1$ and $N > 2$. Let us fix a generator $y$ of $\Gamma$. Let $q$ be a
primitive $N$-th root of 1. We consider the following datum of quantum
linear plane: $$g_{1} = g_{2} = y \in \Gamma, \quad \chi_{1}, \chi_{2} \in
\widehat{\Gamma}, \quad \chi_{1}(y) = q,\quad \chi_{2}(y) = q^{-1}. $$
We consider the compatible datum $$\Cal D = \left(\mu_{1} = 1,
\quad \mu_{2} = 1, \quad \lambda_{ij} =  \lambda\right),$$
where $\lambda\in k$ is arbitrary. 

As above, given a positive integer $n$, $\Bbb G_{n}$ denotes the group of
$n$-th roots of 1 in $k$.

\proclaim{Theorem  6.1} Let $\bq$ be the algebra presented by generators
$h$, $a_{1}$, $a_{2}$ with defining relations

\roster \item"(6.1)"  $h^{NM} = 1$; 

\item"(6.2)" $ha_{1} = qa_{1}h$, $ha_{2} = q^{-1}a_{2}h$;

\item"(6.3)" $a_{1}^{N} = 1 - h^{N}$, $a_{2}^{N} = 1 - h^{N}$;

\item"(6.4)" $a_{2}a_{1} - q a_{1}a_{2} = \lambda \left(1 - h^{2} \right)$.
\endroster 
Then $\bq$ has dimension $MN^{3}$ and carries a Hopf algebra structure given
by
$$\Delta(h) = h\otimes h,
\qquad \Delta(a_{i}) = a_{i} \otimes 1 + h\otimes a_{i}, \quad 1\leq i \leq
2. $$ 
It is pointed and its coradical filtration is given by

\medpagebreak
$$\bq_{n} = \langle h^{i} a_{1}^{j_{1}}
a_{2}^{j_{2}}:\; 0 \leq i \leq NM,
\;   0 \leq j_{1},\;  0 \leq j_{2}, \;  j_{1} + j_{2} \leq n\rangle.
\tag 6.5$$ In particular, 
$$\aligned 
P_{h, 1} (\bq) &= k(1-h)\oplus ka_{1} \oplus ka_{2} \\
P_{g, 1} (\bq) &= k(1-g) \quad \text{ if } g\in \Gamma, \; g\neq h. 
\endaligned \tag 6.6$$

The Hopf algebras $\bq$ and $\bqt$ are isomorphic if and only if
$\widetilde{\lambda} = u\lambda$ for some $u \in \Bbb G_{N}$.
\endproclaim
\demo{Proof} The Hopf algebra structure and the dimension statements
follow from Lemma 5.1 and Proposition 5.2. The description of the coradical
follows from Corollary 5.3. 

We prove now the isomorphism statement. We denote by $\widetilde{h}$,
$\widetilde{a_{i}}$, the generators of $\bqt$. We assume first that $\bq$
and $\bqt$ are isomorphic; let $\phi: \bq \to \bqt$ be a Hopf algebra
isomorphism. Then $\phi$ induces a linear isomorphism 
$$P_{h, 1} (\bq) @>\sim>> P_{\phi(h), 1} (\bqt).$$
By (6.6), $\dim P_{h, 1} (\bq) = 3$; hence $\dim P_{\phi(h), 1} (\bqt)
= 3$ and by (6.6) again, we have 
$\phi(h) = \widetilde{h}.$

Let us write $\phi(a_{1}) = \alpha_{1} (1- \widetilde{h}) + \alpha_{2}
\widetilde{a_{1}} + \alpha_{3} \widetilde{a_{2}}$, for some $\alpha_{i} \in
k$. 

By (6.2), we have $\phi(h) \phi(a_{1}) \phi(h)^{-1} = q \phi(a_{1})$.  Hence
$\alpha_{1} = 0 = \alpha_{3}$ and  $\phi(a_{1}) = \alpha_{2}
\widetilde{a_{1}}$, with $\alpha_{2} \neq 0$. By a similar reason,
$\phi(a_{2}) = \beta_{3} \widetilde{a_{2}}$, with $\beta_{3} \neq 0$. 
Now, by (6.3), $$
1 - \widetilde{h}^{N} = \phi(1 - h^{N}) = \phi(a_{1}^{N})
=\alpha_{2}^{N} \widetilde{a_{1}}^{N} =  \alpha_{2}^{N} \left(1 -
\widetilde{h}^{N}\right).$$
Hence $\alpha_{2}^{N} = 1$, and similarly $\beta_{3}^{N} = 1$. Notice finally
that (6.4) implies
$$\alpha_{2}\beta_{3}\widetilde{\lambda} = \lambda.$$ 
Conversely suppose that $\widetilde{\lambda} = u\lambda$ for some $u \in
\Bbb G_{N}$. Then there is a Hopf algebra isomorphism $\phi: \bq \to \bqt$
uniquely determined by
$$\phi(h) = \widetilde{h}, \quad \phi(a_{1}) =  \widetilde{a_{1}}, \quad
\phi(a_{2}) =  u\widetilde{a_{2}}. \qed$$\enddemo

The following result is a consequence of the argument of the proof of the
Theorem and answers a question of A. Masuoka. 

\proclaim{Corollary 6.2}The group of Hopf algebra automorphisms of $\bq$
is finite. \endproclaim
\demo{Proof} Indeed, any automorphism $T$ has the following form, for
some $j \in \Bbb Z/N$:
$$
T(h) = h, \qquad T(a_{1}) = q^{j}a_{1}, \qquad T(a_{2}) = q^{-j}a_{2}.
\qed$$ \enddemo

\remark{Remark}The Hopf algebra $\bq$ arises as a central extension: $$1
@>>>k[h^{N}] @>>> \bq @>\pi>> {\Cal A}\left(\widehat{\Gamma},
\widehat{\Cal R}, \widehat{\Cal D}\right)  @>>> 1,$$
but $\pi$ has no Hopf algebra section. As $M$ and $N$ could be coprime, this 
shows that Zassenhaus theorem does not generalize to Hopf algebras.
\endremark

\demo{Proof of Theorem 0.3} It is an immediate consequence of Theorem
6.1, letting $M = N = p$. \qed\enddemo 

\remark{Remark} There are also easy examples with $\Gamma =
\Bbb Z/NM_{1} \oplus \Bbb Z/NM_{2}$ of families of pointed non-isomorphic
Hopf algebras of dimension $N^{4}M_{1}M_{2}$, in particular of dimension
$p^{6}$. The construction and proof are very similar. \endremark

\bigpagebreak
\subhead \S 7. Pointed Hopf algebras of order $p^{3}$\endsubhead 

Let $A$ be a non-cosemisimple  pointed Hopf algebra of order
$p^3$, and let $\Gamma$ be the group of its  group-like elements.
By Nichols-Zoeller Theorem \cite{NZ},  we have the following
possibilities: $$ {\text (i)}\; \Gamma = \Bbb Z/(p) \times \Bbb
Z/(p), \qquad {\text (ii)}\; \Gamma = \Bbb Z/(p^2), \qquad {\text
(iii)}\; \Gamma = \Bbb Z/(p).$$
We shall discuss the cases separately and deduce from Theorem 0.2 that in
case (i) $A$ should be of type (a),   in case (ii) $A$ should be of type (b), (c)
or (d) and in case (iii) $A$ should be of type (e) or (f).

\demo{Case (i)} Here $k(\Gamma)$ has index $p$ in $A$ and Theorem 0.2  (ii)
applies. Relation (0.3) turns to $a^{p} = 0$, because any element in
$\Gamma$ has order $p$. It is easy to see that $A \simeq k(\ker
\chi)\otimes k\langle g, a\rangle$, and that the second factor is
isomorphic to a Taft algebra. \qed\enddemo

\bigpagebreak
\demo{Case (ii)} Again, $k(\Gamma)$ has index $p$ in $A$ and Theorem 0.2  (ii)
applies. Let $g, \chi, q, a$  be as in Theorem 0.2  (ii); the order of $q$ is $p$.

We assume first that the order of $g$ is also $p$. Then the relation (0.3)
implies $a^p =0$. On the other  hand, let $h\in \Gamma$ be the generator such
that $h^p = g$. Clearly,  $\xi := \chi(h)$ has order $p^2$. We claim that there
is an isomorphism  of Hopf algebras $$A \simeq k\langle h, x\vert \;
hxh^{-1} = \xi x, \; h^{p^2} = 1, \; x^p = 0\rangle,$$  where the 
comultiplication in the right hand side is as in type (b). Indeed the existence
of a surjective homomorphism from the right hand side to the left  follows
from the considerations above; by a dimension argument it is  an
isomorphism. So, we are in type (b).

\bigpagebreak We assume next that the order of $g$ is $p^2$.
Hence, $a^p = \lambda (1-g^p)$ for some $\lambda \in k$. If
$\lambda = 0$, $A$ is of type (c); otherwise we replace $a$ by
$(^p\sqrt{\lambda})^{-1}a$ and conclude that $A$ is of type
(d). \qed\enddemo

\bigpagebreak 
\demo{Case (iii)} Now $k(\Gamma)$ has index $p^{2}$ in $A$ and Theorem 0.2 
(iii) applies. We observe that possibility (a) is excluded, since every element
of $\Gamma$ has order $p$. Let $g_{i}, \chi_{i}, q_{i}, a_{i}$  be as in
Theorem 0.2  (iii). We set $g = g_{1}$ and  $q = \chi(g)\in \Bbb G_p$. There
are integers  $m, n$ such that $g_{2} = g^{m}$ and $\chi_{2}(g) = q^{n}$. But
$\chi_{1} (g_{2}) \chi_{2} (g_{1}) = 1$ forces $n = -m$.

\bigpagebreak
Relations (0.7) turn to $a_{i}^{p} = 0$. If $\lambda = 0$ in (0.9), then $A$ is
isomorphic to a book algebra and  is of type (f).  If $\lambda \neq 0$, then
$\chi_{1}\chi_{2} = 1$ implies $m = -1$. It is now clear that $A$ is
isomorphic to the Frobenius-Lusztig kernel; that is, it  is of type (e). 
\qed\enddemo


\newpage



\Refs\widestnumber\key{AAA1}

\ref \key  A  \by N. Andruskiewitsch  \paper Notes on extensions
of Hopf algebras \jour Canad. J. Math.   \vol  48  \yr  1996 \pages
3--42 \endref

\ref \key  AS1  \by N. Andruskiewitsch and H.-J. Schneider\paper
Simplicity of some pointed Hopf algebras \paperinfo Appendix to
\cite{A} \endref

\ref \key  AS2  \bysame \paper Hopf Algebras of order $ p^2$ and
braided Hopf algebras of order $p$\jour J. Algebra\toappear
\endref

\ref \key  AS3  \bysame \paper On Hopf Algebras whose
coradical is a Hopf subalgebra\paperinfo in preparation \endref

\ref \key  Be  \by G. Bergman  \paper The Diamond Lemma for Ring Theory 
\jour Adv.  Math.   \vol  29  \yr  1978 \pages 178--218 \endref

\ref \key B \by  N. Bourbaki \book Alg\'ebre commutative.
Chapitre III   \publ Hermann \yr  1961 \endref

\ref \key CM \by W. Chin and I. Musson \paper The coradical
filtration for quantized universal enveloping algebras \jour J.
London Math. Soc.\vol 53\yr 1996\issue 2\pages 50--67 \endref

\ref \key K \by  I. Kaplansky \book Bialgebras \publ University of
Chicago\yr 1975 \endref




\ref \key LR \by R.G. Larson and D.E. Radford \paper  Finite
dimensional cosemisimple Hopf algebras in  characteristic 0 are
semisimple\jour  J. Algebra\vol 117 \yr 1988\pages 267--289
\endref


\ref \key L1 \by G. Lusztig\paper Finite dimensional Hopf
algebras arising from quantized universal enveloping algebras
\jour J. of Amer. Math. Soc.\vol 3\issue 1\pages 257--296 \endref
 
\ref \key L2 \bysame \paper Quantum groups  at  roots  of 1 \jour
Geom. Dedicata \yr 1990 \vol 35 \pages 89--114\endref

\ref \key Mj \by S. Majid  \paper  Crossed products by braided
groups and bosonization \jour J. Algebra\vol 163  \yr 1994\pages
165--190 \endref

\ref \key Ma \by A. Masuoka  \paper  Self dual Hopf algebras of
dimension $p^3$ obtained by extension \jour J. Algebra\vol 178 
\yr 1995\pages 791--806 \endref

\ref \key M \by  S. Montgomery \book  Hopf algebras and their
actions on rings  \publ AMS \yr  1993 \endref


\ref \key N\by W.D. Nichols   \paper Bialgebras of type one  \jour
Commun.  Alg. \yr 1978 \vol 6\pages 1521--1552 \endref



\ref \key NZ \by W.D. Nichols and   M.B. Zoeller   \paper A Hopf
algebra freeness theorem  \jour Amer.J.  of Math. \yr 1989 \vol
111\pages 381--385
 \endref

\ref \key OS \by U. Oberst and H.-J. Schneider  \paper \"Uber
Untergruppen endlicher algebraischer Gruppen   \jour 
Manuscripta Math. \vol 8  \yr 1973 \pages 217--241 \endref



\ref \key R1 \by D. Radford  \paper On the coradical of a finite
dimensional
 Hopf algebra
 \jour Proc.  Amer. Math. Soc.\vol 53 \yr 1975 \pages 9--15 \endref

\ref \key R2 \bysame  \paper  The order of  the antipode of a
finite dimensional Hopf algebra is finite
 \jour Amer.J.  of Math.\vol 98  \yr 1976 \pages 333--355 \endref

\ref \key R3 \bysame  \paper  Hopf algebras with projection
 \jour J. Algebra\vol 92  \yr 1985 \pages 322--347 \endref

\ref \key  S  \by  H.-J. Schneider\paper Finiteness results for semisimple
Hopf algebras \paperinfo manuscript \endref

\ref \key St \by D. Stefan\paper The set of types of n-dimensional
semisimple and cosemisimple   Hopf algebras is finite
 \jour J. Algebra\vol 193  \yr 1997 \pages 571--580 \endref

\ref \key Sw \by Sweedler, M.  \book Hopf algebras
 \yr 1969\publ Benjamin \publaddr New York \endref

\ref \key TW \by E. Taft and  R. L. Wilson\paper On antipodes in
pointed   Hopf algebras 
 \jour J. Algebra\vol 29  \yr 1974 \pages 27--32 \endref

\ref \key T \by M. Takeuchi  \paper  Some topics on $GL_q(n)$
\jour J. Algebra \vol 147  \yr 1992\pages 379--410 \endref


\ref \key Z \by Y. Zhu  \paper Hopf algebras  of prime dimension 
\jour Internat. Math. Res. Notes \yr 1994 \vol 1 \pages 53--59
\endref

 \endRefs
\enddocument